\documentclass[11pt]{amsart}

\usepackage{amsmath} \usepackage{graphicx}

\newcommand{\alg}{\ensuremath{\mathcal{A}}}
\newcommand{\ms}{\ensuremath{\mathcal{M}}}
\newcommand{\df}{\ensuremath{\partial}}
\newcommand{\dbar}{\ensuremath{\bar{\partial}}}

\newcommand{\sgn}{\mbox{sgn}~}

\newcommand{\cc}{\ensuremath{\mathbb{C}}}
\newcommand{\rr}{\ensuremath{\mathbb{R}}}
\newcommand{\zz}{\ensuremath{\mathbb{Z}}}
\newcommand{\qq}{\ensuremath{\mathbb{Q}}}
\newcommand{\pp}{\ensuremath{\mathbb{P}}}

\newcommand{\mylabel}[1]{\label{#1}
}

\theoremstyle{plain}
\newtheorem{theorem}{Theorem}[section]
\newtheorem{corollary}[theorem]{Corollary}
\newtheorem{lemma}[theorem]{Lemma}

\newtheorem{proposition}[theorem]{Proposition}

\theoremstyle{definition}
\newtheorem{definition}[theorem]{Definition}

\theoremstyle{remark}
\newtheorem*{remark}{Remark}

\numberwithin{equation}{section}

\begin{document}

\title{Invariants of Legendrian Knots and  Coherent Orientations}

\author[J. Etnyre]{John B. Etnyre} \address{Stanford University,
  Stanford, CA 94305} \email{etnyre@math.stanford.edu}
\urladdr{http://math.stanford.edu/\textasciitilde etnyre} \thanks{JBE
  is partially supported by an NSF Postdoctoral Fellowship (Grant \#
  DMS-0072853)}

\author[L. Ng]{Lenhard L. Ng} \address{Massachusetts Institute of
  Technology, Cambridge, MA 02139} \email{lenny@math.mit.edu}
\urladdr{http://math.mit.edu/\textasciitilde lenny} \thanks{LLN is
  partially supported by grants from the NSF and DOE}

\author[J. Sabloff]{Joshua M. Sabloff}
\address{Stanford University, Stanford, CA 94305}
\email{sabloff@math.stanford.edu}
\urladdr{http://math.stanford.edu/\textasciitilde sabloff}
\thanks{JMS is partially supported by an NSF Graduate
  Student Fellowship and an ARCS Fellowship}



\begin{abstract}
  We provide a translation between Chekanov's combinatorial theory for
  invariants of Legendrian knots in the standard contact $\rr^3$ and a
  relative version of Eliashberg and Hofer's contact homology.  We use
  this translation to transport the idea of ``coherent orientations''
  from the contact homology world to Chekanov's combinatorial setting.
  As a result, we obtain a lifting of Chekanov's differential graded
  algebra invariant to an algebra over $\zz[t, t^{-1}]$ with a full
  \zz\ grading.
\end{abstract}

\bibliographystyle{plain}
\maketitle

\section{Introduction}
\mylabel{sec:intro}

Legendrian and transversal knot theory has had an extensive influence
on the study of contact $3$-manifolds.  Early on, Bennequin
\cite{bennequin} discovered exotic ``overtwisted'' contact structures
on $\rr^3$ using transversal knots.  Later authors have used
Legendrian and transversal knots to detect overtwisted structures and
to construct and distinguish tight ones (see
\cite{etnyre-honda:non-existence, etnyre-honda:non-filling, kanda:T3,
  lisca-matic}, for example).  On the topological side, work by
Rudolph \cite{rudolph} and others (see \cite{kron-mrowka:contact,
  lisca-matic}) has linked invariants of Legendrian knots of a given
knot type with its smooth slicing properties.

In this paper, we will develop tools for distinguishing Legendrian
knots in the standard contact $\rr^3$.  A natural question to ask is
whether Legendrian knots in a fixed oriented smooth knot type are
classified by their Thurston-Bennequin invariant, $tb,$ and rotation
number, $r.$ If so, then we call that smooth knot type
\textbf{Legendrian simple}.  Though it was never widely believed that
all knot types were Legendrian simple, early evidence and a lack of
suitable invariants suggested that this might be the case. By studying
characteristic foliations on spanning disks, Eliashberg and Fraser
\cite{yasha-fraser} proved unknots are Legendrian simple in the early
1990s.  A few years later, Fuchs and Tabachnikov \cite{f-t} proved
that $tb$ and $r$ are the only finite-type invariants of Legendrian
knots.  More recently, Etnyre and Honda \cite{etnyre-honda:knots}
proved that torus knots and the figure eight knot are also Legendrian
simple.

In the mid-1990's, Chekanov \cite{chv} developed a method of
associating a differential graded algebra (\textbf{DGA}) over $\zz/2$
to the $xy$ diagram of a Legendrian knot $K$.  The generators of this
DGA correspond to the crossings of the diagram and the differential
comes from counting certain immersed polygons whose edges lie in the
diagram of the knot and whose vertices lie at the crossings.  He
proved, combinatorially, that the ``stable tame isomorphism class'' of
the DGA (see Section~\ref{ssec:algebraic} for the precise definition)
is invariant under Legendrian isotopy.  He then proceeded to find an
example of two Legendrian realizations of the $5_2$ knot that have the
same $tb$ and $r$, yet are not Legendrian isotopic.

The first goal of this paper is to lift Chekanov's DGA from $\zz/2$ to
$\zz[t, t^{-1}]$ coefficients and to provide the DGA with a \zz\
grading, regardless of the rotation number of the knot.
Chekanov's original
DGA can be recovered by setting $t=1,$ which will force the grading to
be reduced modulo $2r,$ and taking the coefficients modulo $2$.
Part I of the paper is devoted to this goal.

Concurrent with Chekanov's work on his DGA, Eliashberg and Hofer
adapted the ideas of Floer homology to the contact setting.  Though we
will flesh out a relative version of their ``contact homology theory''
in Section~\ref{sec:ch}, the story goes roughly as follows: let $(M,
\alpha)$ be a contact manifold with a Legendrian submanifold $K$. Let
\alg\ be the free associative unital algebra generated by the Reeb
chords --- i.e. Reeb trajectories that begin and end on $K$.  The
generators are graded by something akin to the Maslov index.  There is
a differential on \alg\ that comes from counting rigid $J$-holomorphic
disks in the \textbf{symplectization} $(M \times \rr, d(e^\tau
\alpha))$ of $M$.  Here, $J$ is a vertically-invariant almost complex
structure compatible with $d(e^\tau\alpha).$ Using Floer and Hofer's
idea of coherent orientations \cite{egh, floer-hofer}, it is possible
to orient all of the moduli spaces of rigid $J$-holomorphic disks used
in the definition of the differential.  As a result, we may use \zz\
coefficients in the definition of the algebra \alg.

The second goal of this paper, carried out in Part II, is to prove
that Chekanov's DGA, and our generalization of it, is a combinatorial
translation of relative contact homology.  Knowing the relation
between the combinatorial and geometric versions of contact homology
is quite useful.  In particular, the lifting of Chekanov's DGA from
$\zz/2$ to $\zz[t,t^{-1}]$ was accomplished by studying this
relationship. Moreover, explicit computations in the framework of
Eliashberg and Hofer's contact homology theory can be difficult, while
computations in the combinatorial theory are more straightforward.
Thus, our translation between the two theories yields many explicit
computations in contact homology.

The paper consists of essentially two parts.  After recalling several
basic ideas from contact geometry in Section~\ref{sec:basics}, we
proceed, in Part I, to describe the combinatorial theory. This part is
self-contained apart from a few technical proofs that are relegated to
an appendix of Part I.  In Part II of the paper we discuss Eliashberg
and Hofer's contact homology and coherent orientations. We then prove
the combinatorial theory developed in Part I is a faithful translation
of this more geometric theory.

\section{Basic Notions}
\mylabel{sec:basics}

We begin by describing some basic notions in three-dimensional contact
geometry.  A \textbf{contact structure} on a $3$-manifold $M$ is a
completely non-integrable $2$-plane field $\xi$.  Locally, a contact
structure is the kernel of a $1$-form $\alpha$ that satisfies the
following non-degeneracy condition at every point in $M$: $$\alpha
\wedge d\alpha \neq 0.$$
In this paper, we will be interested in the
standard contact structure $\xi_0$ on $\rr^3$, which is defined to be
the kernel of the $1$-form $$\alpha_0 = dz + x\,dy.$$
To each contact
form $\alpha$, we may associate a \textbf{Reeb field} $X_\alpha$ that
satisfies $d\alpha(X_\alpha, \cdot)=0$ and $\alpha(X_\alpha)=1$.  By
Darboux's theorem, every contact manifold is locally contactomorphic
to $(\rr^3, \xi_0)$.  See \cite[chapter 8]{aeb} for an introduction to
the fundamentals of contact geometry.

Our primary objects of study are \textbf{Legendrian knots} in $\rr^3$,
i.e. knots that are everywhere tangent to the standard contact
structure $\xi_0$.  In particular, we examine Legendrian isotopy
classes of Legendrian knots, in which two knots are deemed equivalent
if they are related by an isotopy through Legendrian
knots.

Legendrian knots are plentiful: it is not hard to prove that any
smooth knot can be continuously approximated by a Legendrian knot.
Put another way, every smooth knot type has a Legendrian
representative.  The interactions between Legendrian and smooth
knot types constitute a rich and subtle subject.  The first step
in analyzing this interaction is to introduce the ``classical''
invariants $tb$ and $r$ for Legendrian knots in $(\rr^3, \xi_0)$.
(See \cite[chapter 8]{aeb} for more general
  definitions than we give here.)  The \textbf{Thurston-Bennequin}
invariant measures the self-linking of a Legendrian knot $K$.  More
precisely, let $K'$ be a knot that has been pushed off of $K$ in a
direction tangent to the contact structure.  Define $tb(K)$ to be the
linking number of $K$ and $K'$.  The \textbf{rotation number} $r$ of
an oriented Legendrian knot $K$ is the rotation of its tangent vector
field with respect to any global trivialization of $\xi_0$ (e.g.
$\{\partial_x, \partial_y - x\,\partial_z \}$).

In this paper, we use the combinatorics of generic projections of
Legendrian knots into the $xy$ plane extensively.  Not all knot
diagrams in the $xy$ plane can be lifted to Legendrian knots: Stokes'
Theorem implies that the diagram must bound zero (signed) area.
Chekanov describes the combinatorial restrictions on the form of the
$xy$ projection of a Legendrian knot in \cite{chv}. Note that the
Thurston-Bennequin number of $K$ may be computed from the writhe of
the $xy$ projection of $K$ while the rotation number of $K$ is just
the (counterclockwise) rotation number of the diagram.  For example,
the Legendrian unknot in Figure~\ref{fig:word} has $tb = -2$ and
$r=1$.

\bigskip
\begin{center}
{\bf
Part I \break
The Combinatorial Theory
}
\end{center}
\medskip

\noindent
In Section~\ref{sec:alg} of the paper we describe the contact homology
of a knot in $\rr^3$ in purely combinatorial terms. We do this by
giving a self-contained generalization of Chekanov's differential
graded algebra. In broad outline we follow \cite{chv}, making the
necessary changes to extend the DGA defined there over $\zz/2$ to a
DGA over $\zz[t,t^{-1}].$ We provide several illustrative computations
in Section~\ref{sec:ex}.

\section{A Combinatorial Definition of the Algebra}
\mylabel{sec:alg}

Given an oriented Legendrian knot $K$ in standard contact structure on
$\rr^3$ we show how to associate a differential graded algebra over
$\zz[t,t^{-1}]$ to it in Sections~\ref{ssec:alg} and \ref{ssec:d}.
Then, in Section~\ref{ssec:algebraic}, we describe an equivalence
relation, stable tame isomorphism, on DGA's and show that the
equivalence class of the DGA associated to a Legendrian knot is
invariant under Legendrian isotopy. This in turn implies that the
homology of the DGA is invariant under Legendrian isotopy. The
algebras and homology that we work with are non-abelian and hence
somewhat hard to use.  In Section~\ref{ssec:abelian} we prove that
when the DGA of a Legendrian knot is abelianized (over $\qq$), its
homology is still invariant under Legendrian isotopy.

\subsection{From the Knot $K$ to the Algebra \alg}
\mylabel{ssec:alg}

We begin by decorating a generic $xy$ diagram of a given Legendrian
knot $K$.  First, label the crossings of $K$ by $\{a_1, \ldots,
a_n\}$.  Next, label each
quadrant around a crossing as shown in Figure~\ref{fig:pq}.
\begin{figure}
  \centerline{\includegraphics{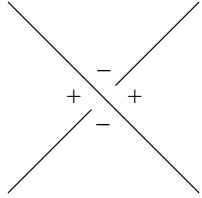}}
  \caption{The Reeb signs around a crossing.}
  \mylabel{fig:pq}
\end{figure}
We will
refer to these labels as the \textbf{Reeb signs}.
Call a quadrant at a crossing \textbf{positive} or \textbf{negative}
depending on its Reeb sign.

\begin{definition} \mylabel{defn:alg}
  The algebra $\alg(a_1, \ldots, a_n)$ is the graded free associative
  unital algebra over $\zz[t, t^{-1}]$ generated (as an algebra) by
  $\{a_1, \ldots, a_n\}$.
\end{definition}
The grading for $t$ is defined to be $2r(K)$.  To grade a generator
$a_i$, we first need to specify a capping path $\gamma_{a_i}$.
By this we mean  the unique path $\gamma_{a_i}$ in $K$ which begins at
the undercrossing of $a_i$,
follows in the direction of the orientation of $K$, and ends when it
reaches the overcrossing of $a_i$; see Figure~\ref{fig:cw-ccw}.
\begin{figure}
  \centerline{\includegraphics{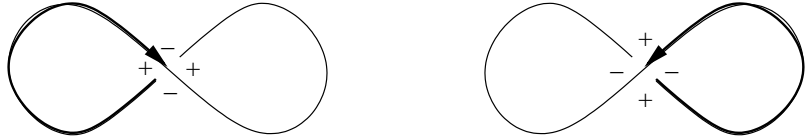}}
\vspace{12pt}
  \caption{The choice of capping path $\gamma$ for a crossing.
The path $\gamma$ is denoted by a heavy line; the arrows indicate
the orientation of the knot and of $\gamma$; and the signs
are Reeb signs.  The diagrams show a crossing
coherent about a $+$ and about a $-$, respectively.
}
  \mylabel{fig:cw-ccw}
\end{figure}
Let the rotation number $r(\gamma_{a_i})$ be the fractional
number of counterclockwise revolutions made by a
tangent vector to $\gamma_{a_i}$
as we traverse the path.  (Thus the rotation number of the entire oriented
knot $K$ is the usual rotation number.)
We may perturb the diagram of $K$ so that
all crossings are orthogonal; then $r(\gamma_{a_i})$ is an odd multiple of
$1/4$.  Define
\begin{equation}
  \mylabel{eqn:grading}
  |a_i| = -2r(\gamma_{a_i}) - \frac{1}{2}.
\end{equation}
In addition, define the \textbf{sign of a crossing} $a_i$ to be $\sgn
a_i = (-1)^{|a_i|}$.  (Note that we may recover Chekanov's original
grading by setting $t=1$, but this forces us to consider the grading
modulo $2r(K)$.)

We need one more piece of notation.  Define $a_i$ to be
\textbf{coherent about a $+$} (respectively \textbf{coherent about a
  $-$}) if, in a neighborhood of $a_i$, the quadrant enclosed by the
path $\gamma_{a_i}$ is labeled by a $+$ (resp.\ $-$).  We also say
that the quadrants that are enclosed by $\gamma_{a_i}$ or its
complement are \textbf{coherent}; the remaining two quadrants are
\textbf{incoherent}.  The following lemma will be useful in the
future.

\begin{lemma}
  \mylabel{lem:coherent-grading} If a crossing is coherent about a
  $+$ (resp.\ $-$), then the grading of the associated variable
  $a_i$ is odd (resp.\ even).
\end{lemma}

\begin{proof}
  For a crossing coherent about a $+$, the rotation
  number of $\gamma$ is $k - \frac{3}{4}$ for some integer $k$ (see Figure~\ref{fig:cw-ccw}); then by
  definition,
  \begin{equation}
    \mylabel{eqn:parity}
    |a_i| = -2(k-\frac{3}{4}) - \frac{1}{2} = -2k+1.
  \end{equation}
  The proof for the other case is identical.
\end{proof}

\subsection{The Differential $\df$}
\mylabel{ssec:d}

We describe the differential $\df$ by appropriately counting immersed
disks in the $xy$-projection of $K.$ For this, we need some notation.
Let $D^2_*=D^2\setminus\{x, y_1,\ldots, y_n\}$, where $\{x,
y_1,\ldots, y_n\}\subset \partial D^2$ are called boundary punctures.
Throughout this section, $\rr^2$ refers to the $xy$-plane.
\begin{definition}
  \mylabel{defn:delta-combinatorial} Fix a homology class $A\in
  H_1(K)=\zz$ and define $\Delta^A(a;b_1, \ldots, b_n)$ to be the
  space of all orientation-preserving immersions $f:(D^2_*,\partial
  D^2_*) \to (\rr^2, K)$ (up to reparametrization) that satisfy:
\begin{enumerate}
\item The homology class $\left[ \left( \pi_{xy}|_K \right) ^{-1}
    \left( \mathrm{Im} (f|_{\partial D^2}) \cup \gamma_a \cup
      -\gamma_{b_1} \cup \cdots \cup -\gamma_{b_n} \right) \right]$
  coincides with $A$.
\item The map $f$ sends the boundary punctures to the crossings of the
  diagram of $K,$ and at a boundary puncture, the map $f$ covers
  either one or three quadrants, with the majority of the quadrants
  positive at the crossing $a$ and negative at the crossings $b_i$.
\end{enumerate}
Formally, we define the dimension of $\Delta^A(a;b_1, \ldots, b_n)$ to be
\begin{equation}\label{eq:dimension}
        \dim (\Delta^A(a; b_1, \ldots, b_n))=|a| - \sum |b_i| + 2n(A)r(K) - 1,
\end{equation}
where $n(A)\in\zz$ is the image of $A$ under the isomorphism $H_1(K)\cong\zz$ given
by the choice of orientation for $K.$
\end{definition}

We call a boundary puncture of $f$ a {\bf convex} or {\bf
non-convex corner}, depending on whether $f$ covers one or three
quadrants, respectively, at the puncture. The formal dimension of
$\Delta^A(a;b_1, \ldots, b_n)$ dictates the number of non-convex
corners of elements of $\Delta^A(a;b_1, \ldots, b_n)$.

\begin{lemma}
Any $f \in \Delta^A(a;b_1,\ldots,b_n)$ has precisely $\dim (\Delta^A(a; b_1, \ldots, b_n))$
non-convex corners.
\mylabel{lem:nonconvexcorners}
\end{lemma}

\begin{proof}
Suppose that $f$ has $k$ non-convex corners.
Let $\gamma$ be the closed curve in the diagram of $K$ which is the union of the oriented
boundary $f(\df D^2)$ and the paths $\gamma_a,-\gamma_{b_1}, \ldots,-\gamma_{b_n}$.  Then
the rotation number of $\gamma$ is, by definition, $n(A)r(K)$.

  On the other hand, the rotation number of $\gamma$ is also the sum
  of the rotation numbers of its pieces.  We may assume that the
  crossings of the diagram of $K$ are orthogonal. The sum of the
  rotation numbers of the smooth pieces of $\gamma$ is simply $1-(n+1-2k)/4$,
  because $\gamma$ is traversed counterclockwise, and each corner
  contributes a rotation of $1/4$ if convex, and $-1/4$ if non-convex.
  Also, the rotation number of
  $\gamma_a$ is $-(2|a|+1)/4$, and similarly for $\gamma_{b_i}$.  Thus
  the rotation number of $\gamma$ is
  \begin{equation}
    n(A)r(K) = 1 - \frac{n+1-2k}{4} - \frac{2|a|+1}{4} + \sum_{i=1}^n
    \frac{2|b_i|+1}{4} = \frac{-|a|+\sum |b_i| +k + 1}{2}.
  \end{equation}
  This implies that $k = \dim (\Delta^A(a; b_1, \ldots, b_n))$, as desired.
\end{proof}

We can assign a word in \alg\ to each immersed disk as follows:
starting with the first corner after the one covering the $+$
quadrant, the word is a list of the crossing labels of all subsequent
negative corners encountered while traversing the boundary of the immersed
polygon counter-clockwise.  We also associate a sign
to each immersed disk as follows:
\begin{definition}
  \mylabel{def:signs} To each quadrant $Q$ in the neighborhood of a
  crossing $a$, we associate a sign $\varepsilon_{Q,a}$, called the
  \textbf{orientation sign}, determined from
  Figure~\ref{fig:lenny-signs}.  For an immersed disk with one
  positive corner $a$ (with respect to the Reeb signs) and negative
  corners $b_1,\ldots,b_n$, define the orientation sign
  $\varepsilon(a;b_1,\ldots,b_n)$ to be the product of the orientation
  signs over all corners of the disk.
\end{definition}

\begin{figure}
  \centerline{\includegraphics{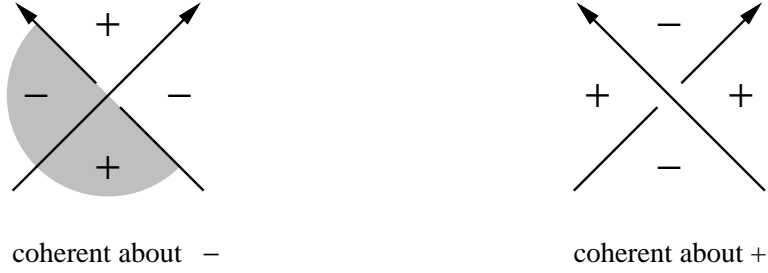}} \vspace{12pt}
  \caption{The signs in the figures are Reeb signs.  The orientation
    signs are $-1$ for
    the two shaded quadrants and $+1$ elsewhere.  }
  \mylabel{fig:lenny-signs}
\end{figure}

Inspection of Figure~\ref{fig:lenny-signs} yields the following
lemma, which will be useful in Sections~\ref{sec:d2} and
\ref{sec:invariance} when we prove that $\df$ is a differential
and \alg\ is invariant under Legendrian isotopy.

\begin{lemma}
  Around a crossing $a$, the product of the orientation signs of two
  opposite quadrants is $-\sgn a$.  The product of the orientation
  signs of two adjacent quadrants is $1$ or $-\sgn a$, depending on
  whether the Reeb signs of the two quadrants are, in counterclockwise
  order around the crossing, $-$ and $+$ or $+$ and $-$, respectively.
  (These cases correspond to the quadrants being on the same side of
  the undercrossing or overcrossing line, respectively.)
  \mylabel{lem:prodsigns}
\end{lemma}

We are ready to define a differential on \alg.
\begin{theorem}
  \mylabel{thm:dga} The algebra \alg\ is a differential graded algebra
  (DGA) whose differential $\df$ is defined as follows:
  \begin{equation}
    \mylabel{eqn:comb-d}
    \df a = \sum_{\dim (\Delta^A(a; b_1, \ldots, b_n)) = 0}
    \varepsilon(a; b_1, \ldots, b_n) t^{-n(A)}b_1 \cdots b_n.
  \end{equation}
  Extend \df\ to \alg\ via $\df(\zz[t,t^{-1}]) = 0$ and the signed
  Leibniz rule:
  \begin{equation} \mylabel{eqn:leibniz}
    \df (vw) = (\df v) w + (-1)^{|v|} v (\df w).
  \end{equation}
  The differential \df\ has degree $-1$ and satisfies $\df^2=0$.
\end{theorem}
\noindent
We remark that Lemma~\ref{lem:nonconvexcorners} implies that
the sum in equation (\ref{eqn:comb-d}) is over immersions with no non-convex corners.
Note that we use $t^{-n(A)}$ rather than $t^{n(A)}$ in equation (\ref{eqn:comb-d});
this convention simplifies notation slightly in examples.

The fact that \df\ has degree $-1$ follows directly from the dimension
formula, equation (\ref{eq:dimension}).  To prove Theorem~\ref{thm:dga}, it
suffices to check that $\df^2=0$; this is done in Section~\ref{sec:d2}.

%
%

\subsection{Algebraic Definitions}
\mylabel{ssec:algebraic}

The definition of \alg\ depends heavily on the choice of the projection of
$K$.  The following definitions, due essentially to Chekanov
\cite{chv}, give a notion of equivalence which reflects the possible
changes in \alg\ due to changes in the knot projection under
Legendrian isotopy.

The first definition picks out a particularly simple set of
isomorphisms of \alg.

\begin{definition}
  \mylabel{defn:elem}
  A graded chain isomorphism
  $$\phi: \alg(a_1, ..., a_n) \longrightarrow \alg(b_1,\ldots, b_n)$$
  is \textbf{elementary} if there is some $j \in \{1, \ldots, n\}$
  such that
  \begin{equation} \mylabel{eqn:elem-isom}
    \phi(a_i) =
    \begin{cases}
      b_i, & i \neq j \\
      \pm b_j + u, & u \in
      \alg(b_1,\ldots,b_{j-1},b_{j+1},\ldots,b_n),\ i=j.
    \end{cases}
  \end{equation}
  A composition of elementary isomorphisms is called \textbf{tame}.
\end{definition}

\noindent We note that since it is not known whether or not all
isomorphisms are tame we should technically be working with
``semi-free'' algebras (free algebras with specified generators).
This technical point will not cause any problems so we ignore it.
In \cite{chv}, this point is carefully discussed and the
interested reader is referred there.

We also need an algebraic operation that reflects the second
Reidemeister move (see move III in Figure~\ref{fig:reidemeister}).
Define a special algebra $\mathcal{E}_i = \alg(e_1, e_2)$ by
setting  $|e_1| = i$, $|e_2| = i-1$, $\df e_1 = e_2$, $\df e_2 =
0$.

\begin{definition} \mylabel{defn:stabilization}
  The degree $i$ \textbf{stabilization} $S_i(\alg(a_1, \ldots, a_n))$
  of $\alg(a_1, \ldots, a_n)$ is defined to be $\alg(a_1, \ldots, a_n,
  e^i_1, e^i_2)$.  The grading and the differential are inherited from
  both \alg\ and $\mathcal{E}_i$.
  Two algebras \alg\ and $\alg'$ are \textbf{stable tame isomorphic}
  if there exist two sequences of stabilizations $S_{i_1}, \ldots,
  S_{i_n}$ and $S_{j_1}, \ldots, S_{j_m}$ and a tame isomorphism
  $$\phi: S_{i_n}( \cdots S_{i_1}(\alg) \cdots ) \longrightarrow
  S_{j_m}( \cdots S_{j_1}(\alg') \cdots ).$$
  Two differential algebras $(\alg; \df)$ and $(\alg'; \df')$ are
  \textbf{stable tame isomorphic} if there is a stable tame isomorphism
  from \alg\ to $\alg'$ that is also a chain map.
\end{definition}

\noindent
This equivalence relation is designed for the following important
theorem.

\begin{theorem}
  \mylabel{thm:invariance} The stable tame isomorphism class of
  $\alg(a_1, \ldots, a_n; \df)$ is invariant under Legendrian isotopy
  of $K$.
\end{theorem}

Chekanov proved this theorem over $\zz/2$; we will prove the $\zz[t,
t^{-1}]$ version of this theorem in Section~\ref{sec:invariance}.  As
a corollary, we obtain a proof that the homology of $\alg(a_1, \ldots,
a_n; \df)$ is an invariant:

\begin{corollary}
  \mylabel{cor:homology} The homology $H(\alg(a_1, \ldots, a_n; \df))$
  is invariant under Legendrian isotopy of $K$.
\end{corollary}

\begin{proof}
  It suffices to prove that homology does not change under
  stabilizations.  Consider the natural inclusion and projection
  $$\alg \stackrel{i}{\longrightarrow} S(\alg)
  \stackrel{\tau}{\longrightarrow} \alg.$$
  On one hand, $\tau \circ i =
  Id_\alg$.  We need to prove that $i \circ \tau$ is chain homotopic
  to $Id_{S(\alg)}$, i.e. that there exists some linear map $H:
  S(\alg) \to S(\alg)$ that satisfies
  \begin{equation}
    \mylabel{eqn:chain-homotopy}
    i \circ \tau - Id_{S(\alg)} = H \circ \df + \df \circ H.
  \end{equation}
  It is not hard to check that the following $H$ satisfies these
  requirements:
  \begin{equation}
    \mylabel{eqn:H}
    H(w) = \begin{cases}
      0 & w \in \alg \\
      0 & w=ae_1b \quad \text{with\ } a \in \alg \\
      (-1)^{|a|+1}ae_1b & w=ae_2b \quad \text{with\ } a \in \alg.
    \end{cases}
  \end{equation}
\end{proof}

\subsection{Abelianization of \alg}
\mylabel{ssec:abelian}

One possible way to simplify calculations with \alg\ is to change the
base ring and abelianize.

\begin{definition}
  Given a Legendrian knot $K$ with crossings labeled $\{a_1, \ldots,
  a_n\}$, let $\tilde{\alg}_\qq(a_1, \ldots, a_n; \df)$ be the free
  graded supercommutative associative unital differential algebra over
  $\qq[t, t^{-1}]$ generated as an algebra by $\{a_1, \ldots, a_n\}$.
  The gradings and the differential \df\ are the same as those defined
  in Section~\ref{ssec:d}.
\end{definition}

\noindent By supercommutative we mean that $wv = (-1)^{|v||w|} vw$.
A key feature of the abelianized algebra over a field not of
characteristic $2$ is that for any generator $a$ of odd degree, $a^2 =
0$.  Note that this is not the case if we abelianize over $\zz/2$.

Just as in the non-abelian case, we have the following results
(cf.\ Theorem~\ref{thm:invariance} and Corollary~\ref{cor:homology}).

\begin{proposition}
  \mylabel{prop:abelian} The stable tame isomorphism class of
  $\tilde{\alg}_\qq$ is an invariant of the Legendrian isotopy class
  of the Legendrian knot $K$.
\end{proposition}

\noindent
The proof of \ref{prop:abelian} is a simple diagram chase.

\begin{theorem}
  \mylabel{thm:ab-inv} The homology of $\tilde{\alg}_\qq$ is an
  invariant of the Legendrian isotopy class of the Legendrian knot
  $K$.
\end{theorem}

\begin{proof}
  The only place where we use the non-commutativity of \alg\ in
  proving invariance is in the definition of the map $H$ in
  Corollary~\ref{cor:homology}.  In the abelianized case over \qq, we
  can redefine $H$ so that it is still a chain homotopy.

  If $i$ is even, then we define $H: S_i(\tilde{\alg}) \to
  S_i(\tilde{\alg})$ as follows:
  \begin{equation} \mylabel{eqn:ab-H-even}
    H(w) = \begin{cases}
      -\frac{1}{k+1}e_1^{k+1}a & w \in e_2 e_1^k \tilde{\alg}(a_1, \ldots, a_n) \\
      0 & \text{otherwise;}
    \end{cases}
  \end{equation}
  if $i$ is odd,
  \begin{equation} \mylabel{eqn:ab-H-odd}
    H(w) = \begin{cases}
      -e_2^{k-1} e_1 \tilde{a} & w \in e_2^k \tilde{\alg}(a_1,
      \ldots, a_n, e_1) \\
      0 & \text{otherwise.}
    \end{cases}
  \end{equation}
  It is easy to check that this $H$ works in the abelianized version of
  the proof of Corollary~\ref{cor:homology}.
\end{proof}

\section{Examples}
\mylabel{sec:ex}

In this section, we compute the DGA for three sample knots with
both orientations.

\subsection{Unknot.}
\mylabel{ssec:unknot}

Consider the oriented unknot in Figure~\ref{fig:word} which
has $tb=-2$ and $r=1.$
The capping paths are given by $\gamma_a = \gamma_1$ and $\gamma_b =
\gamma_4 + \gamma_1 + \gamma_2$.  Assuming orthogonal crossings gives
$r(\gamma_a) = -3/4$ and $r(\gamma_b) = 1/4$; hence $|a| = 1$, $|b|
= -1$, and $|t| = 2$.

\begin{figure}
  \begin{center}
    \centerline{\includegraphics[width=5in]{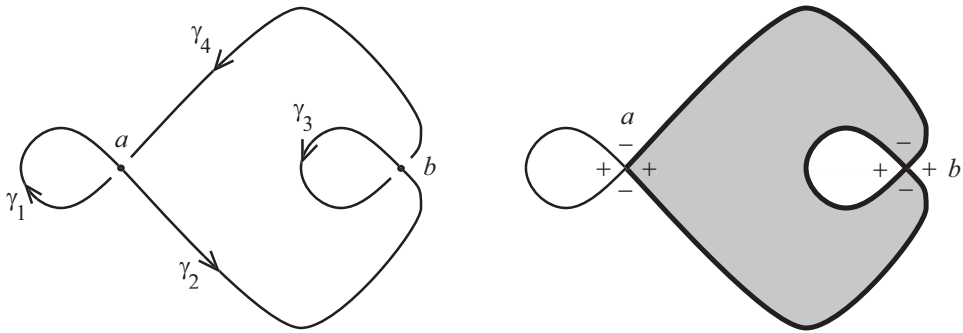}}
\vspace{12pt}
    \caption{The unknot with $tb=-2$ and $r=1.$
On the left, the knot has been divided by crossings $a$ and $b$
into four oriented curves.  On the right, Reeb signs and an embedded
disk are shown.}
    \mylabel{fig:word}
  \end{center}
\end{figure}

The word in $\df a$ represented by the immersed disk shown in
Figure~\ref{fig:word} is
is $t b^2$, where the
power of $t$ follows from the fact that the boundary of the disk is
$\gamma = \gamma_2 - \gamma_3 + \gamma_4$, and
$\gamma + \gamma_a - 2\gamma_b$ winds $-1$ times around the knot.

Continuing in this way, we find that
$\df a = 1 + t b^2$, $\df b = t^{-1}$.  Note that all orientation signs are
positive, because both crossings are odd degree and hence coherent
about a $+$.  Since $t$ has degree $2$, we see that $\df$ does indeed
lower degree by $1$.

If we reverse the orientation of the knot, we may similarly compute that
$r=-1$, $|a| = 3$, $|b| = 1$, $\df a = t^{-1} + b^2$, $\df b = 1$.

This knot is reducible, in the terminology of \cite{chv}; that is,
we can view it as a Legendrian knot with an added loop (the loop around
$b$).  As in \cite{chv}, any reducible knot has a DGA with trivial
homology, and which indeed contains no information besides the
classical invariants.  In our case, for the original orientation,
$\df (tb) = 1$, and so the homology vanishes.

\subsection{Trefoil.}
\mylabel{ssec:trefoil}
Consider the right-handed trefoil knot depicted in
Figure~\ref{fig:examples}.  (This is example 4.3 in \cite{chv}.)
This satisfies $r = 1$, $tb = -6$, and $|a_i| = -1$ for all $i$.
As in the previous example, all orientation signs are positive.

\begin{figure}
  \begin{center}
    \centerline{\includegraphics[width=5in]{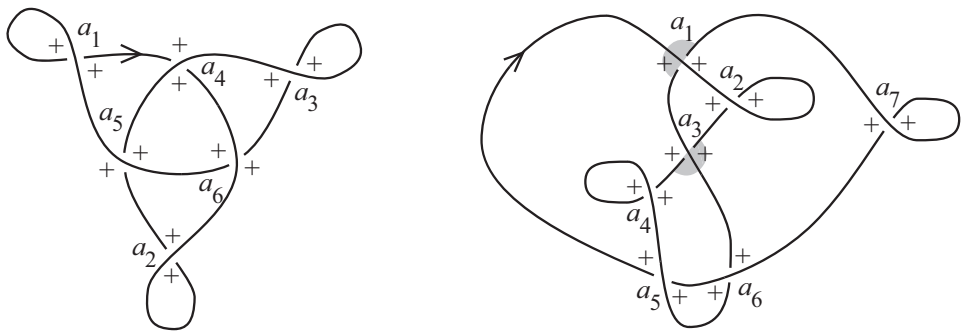}}
\vspace{12pt}
    \caption{Examples of Legendrian trefoil and figure eight knots.
Crossings are labeled, and the $+$ signs represent Reeb signs;
to reduce clutter, the $-$ Reeb signs have been omitted.
The four shaded quadrants are the only quadrants with negative
orientation sign.
}
    \mylabel{fig:examples}
  \end{center}
\end{figure}

We can then compute the differential:
\begin{eqnarray*}
\df a_1 &=& t^{-1} + a_5 a_4 \\
\df a_2 &=& t^{-1} + a_6 a_5 \\
\df a_3 &=& t^{-1} + a_4 a_6 \\
\df a_4 &=& \df a_5 = \df a_6 = 0.
\end{eqnarray*}
Note that, for this knot or any knot with $r \neq 0$, it is easiest
to deduce powers of $t$ from the facts that $t$ has degree $2r$ and
$\df$ lowers degree by $1$.

For the same knot with the opposite orientation, we have
$r = -1$, $|a_i| = 1$ for all $i$, and
\begin{eqnarray*}
\df a_1 &=& 1 + t a_5 a_4 \\
\df a_2 &=& 1 + t a_6 a_5 \\
\df a_3 &=& 1 + t a_4 a_6 \\
\df a_4 &=& \df a_5 = \df a_6 = 0.
\end{eqnarray*}

\subsection{Figure eight.}
\mylabel{ssec:figure8}

Consider the figure eight knot depicted in Figure~\ref{fig:examples}.
We have $r = 0$ and $tb = -3$, and the degrees of the crossings
are given by $|a_2| = |a_4| = |a_5| = |a_7| = 1$,
$|a_1| = |a_3| = 0$, $|a_6| = -1$.

A tip for calculating powers of $t$, especially for knots with $r = 0$:
it is useful to choose
a small section of the knot and count how many times it is traversed
by the boundary of the immersed disk and the capping paths.
For instance, for the figure eight knot, the loop next to $a_7$ is
traversed (positively) by the capping paths for $a_1$, $a_2$, and
$a_6$.  Thus the exponent of $t$ for a monomial in $\df a_i$ is
the number of times $a_1,a_2,a_6$ appears in the monomial, minus one if
$i = 1,2,6$; for the term $1$ in $\df a_7$ only, we must then subtract one,
since this is the only disk whose boundary traverses (positively)
the loop next to $a_7$.

We thus find that
\[
\begin{array}{rclcrcl}
\df a_1 &=& -a_6 + a_6 a_3 + t a_6 a_3 a_5 a_6 & \hspace{0.25in}&
\df a_4 &=& 1 - a_3 - t a_5 a_6 a_3  \\
\df a_2 &=& t^{-1} + a_1 a_3 - a_6 a_3 a_4 &&
\df a_7 &=& t^{-1} + a_3 - t a_3 a_6 a_3 a_5 \\
\multicolumn{7}{c}{\df a_3 = \df a_5 = \df a_6 = 0.}
\end{array}
\]

For the knot with the opposite orientation, the degrees of the crossings
remain the same since $r=0$, but the signs and powers of $t$ in
the differential change:
\[
\begin{array}{rclcrcl}
\df a_1 &=& a_6 + t a_6 a_3 + t a_6 a_3 a_5 a_6 & \hspace{0.25in}&
\df a_4 &=& t^{-1} + a_3 + t a_5 a_6 a_3  \\
\df a_2 &=& 1 + t a_1 a_3 + t^2 a_6 a_3 a_4 &&
\df a_7 &=& 1 - a_3 - t^2 a_3 a_6 a_3 a_5 \\
\multicolumn{7}{c}{\df a_3 = \df a_5 = \df a_6 = 0.}
\end{array}
\]

\bigskip
\begin{center}
{\bf
Appendix to Part I \break
Proofs of $\df^2=0$ and Invariance
}
\end{center}
\medskip

\noindent
Here we show that the arguments in \cite{chv} can be strengthened to prove
Theorems~\ref{thm:dga} and \ref{thm:invariance}.

\section{Proof that $\df^2=0$}
\mylabel{sec:d2}

The geometric motivation behind the following proof will become clear
in Section~\ref{sec:ch}. For now, we give a purely combinatorial proof
that $\df^2 =0.$ Our proof mimics the corresponding proof in
\cite{chv}, so we omit some details and cases but clearly indicate the
complications added by our signs and powers of $t.$

Since \df\ obeys the signed Leibniz rule, it suffices to prove that
$\df^2=0$ on the generators of \alg.  Let $a$ be such a generator.  If
we disregard signs and powers of $t$ for now, a term in $\df^2 a$ is
of the form
\begin{equation} \mylabel{eqn:term1}
  a_1\cdots a_{k} b_2 \cdots b_l c_{1} \cdots c_m,
\end{equation}
where $a_1\cdots a_k b_1 c_1 \cdots c_m$ and $b_2\cdots b_l$ are terms
in $\df a$ and $\df b_1$, respectively.  We may think of this term in
$\df^2 a$ as the two relevant disks glued together at the crossing
$b_1$; see Figure~\ref{fig:obtuse-disk}.
\begin{figure}
  \centerline{\includegraphics[width=5.2in]{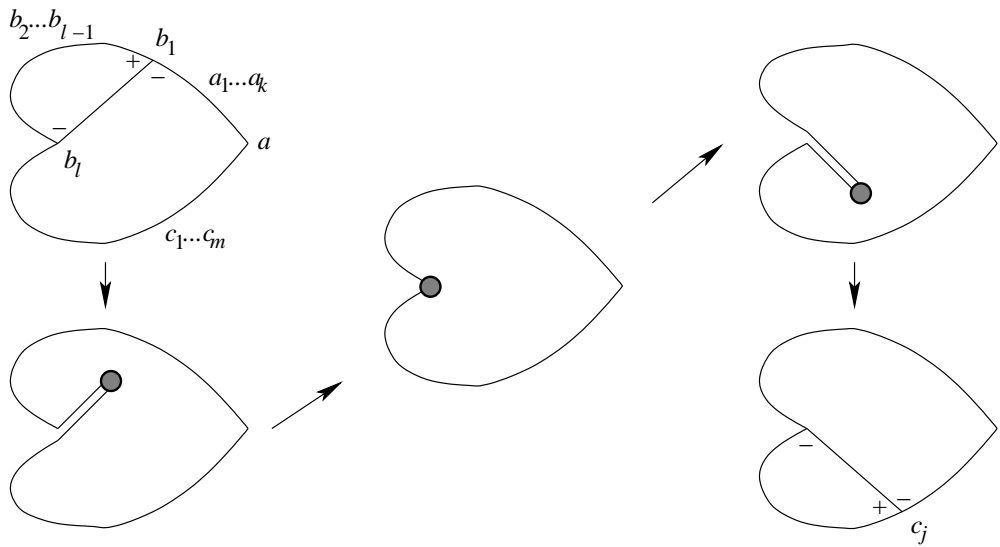}}
  \caption{One possible configuration for disks involved in $\df^2 a$.  The signs are Reeb
    signs.}  \mylabel{fig:obtuse-disk}
\end{figure}
Gluing yields a loop that
has a branch point just after the corner at $b_l$.

After retracting the branch point to the boundary of the obtuse disk
(in Figure~\ref{fig:obtuse-disk}, we retract to $b_l$),
we may think of this loop
as the boundary of an ``obtuse disk,'' i.e., a disk with exactly one
non-convex corner.  Near the non-convex corner, there are two segments
of $K$ pointing into the obtuse disk; we have just pushed the branch
point along one of these segments.  Now push it along the other segment,
and continue until the loop breaks into the boundary of two immersed disks;
see Figure~\ref{fig:obtuse-disk}.
These two immersed disks represent another term contributing to
$\df^2 a$ which cancels our original term, up to signs and powers of $t$.
To complete the proof that $\df^2 = 0$, we need to check that
our two terms share the same power of $t$ and have opposite signs.

We first address the powers of $t$.  As described in Section~\ref{ssec:d},
the (counterclockwise oriented) boundary of each disk, along with the
appropriately oriented capping paths of the corners
of the disk, forms a closed curve in $K$; the negation of the winding
number of this curve around $K$ is the power of $t$ associated to the
disk.  For two disks which glue together to form an obtuse disk,
a quick consideration of the capping paths shows that the powers of
$t$ associated to the two disks multiply to give the power of $t$
associated to the obtuse disk.  Hence the two terms of $\df^2a$
corresponding to a given obtuse disk have the same power of $t$.

It remains to check that the two terms in $\df^2 a$ corresponding to a
fixed obtuse disk have opposite signs.  Consider the obtuse disk $D$
in Figure~\ref{fig:obtuse-disk}, with positive corner at $a$ and
non-convex corner at $b_l$.  For two points $v_1,v_2$ along the
boundary of this disk, let $\sgn \overline{v_1v_2}$ denote the product
of $\sgn w$ over all corners $w$ on the portion of the boundary from
$v_1$ to $v_2$, not including the endpoints; here the boundary is
oriented counterclockwise.

There are several possible configurations, depending on how the paths
from $b_l$ divide $D$.  We will consider one such
configuration, shown in Figure~\ref{fig:obtuse-disk}; the arguments for
the other configurations are similar.

The two terms in $\df^2 a$ arising from $D$ have two sets of
signs: one from the disks themselves through
Definition~\ref{def:signs}, and one from the signed Leibniz rule.
The signs arising from the signed Leibniz rule are $\sgn
\overline{ab_1}$ for the figure at the top left and $\sgn
\overline{ab_l}$ for the figure on the bottom right.  The signs
arising from the disks themselves are the same for both figures,
with the exception of the contribution of the corners marked with
their Reeb signs in Figure~\ref{fig:obtuse-disk}.

Lemma~\ref{lem:prodsigns} applied to Figure~\ref{fig:obtuse-disk}
shows that the total contribution of the marked corners at $b_1$,
$b_l$, and $c_j$ is $-\sgn b_1$, $-\sgn b_l$, and $1$,
respectively. Hence the total sign difference between the two
terms in $\df^2 a$ is
\begin{equation}
  (\sgn \overline{ab_1}) (\sgn \overline{ab_l}) (-\sgn b_1)(-\sgn b_l)
  = (\sgn b_1) (\sgn \overline{b_1b_l}) (\sgn b_l) = -1,
\end{equation}
where the last equality follows from the fact that $\df$ lowers
degree by 1.

This concludes the checking of signs, and the proof that $\df^2 = 0$.
Theorem~\ref{thm:dga} follows. \qed

\section{Invariance of \alg}
\mylabel{sec:invariance}

The goal of this section is to prove Theorem~\ref{thm:invariance}.
Given a diagram of $K$ and its algebra $\alg(a_1, \ldots, a_n)$, we
will check that the stable tame isomorphism class of \alg\ is
invariant under each of the three Legendrian Reidemeister moves (see
Figure~\ref{fig:reidemeister}).  As in Section~\ref{sec:d2}, we will
use an adaptation of Chekanov's original proof over $\zz/2$.

\begin{figure}
   \centerline{\includegraphics{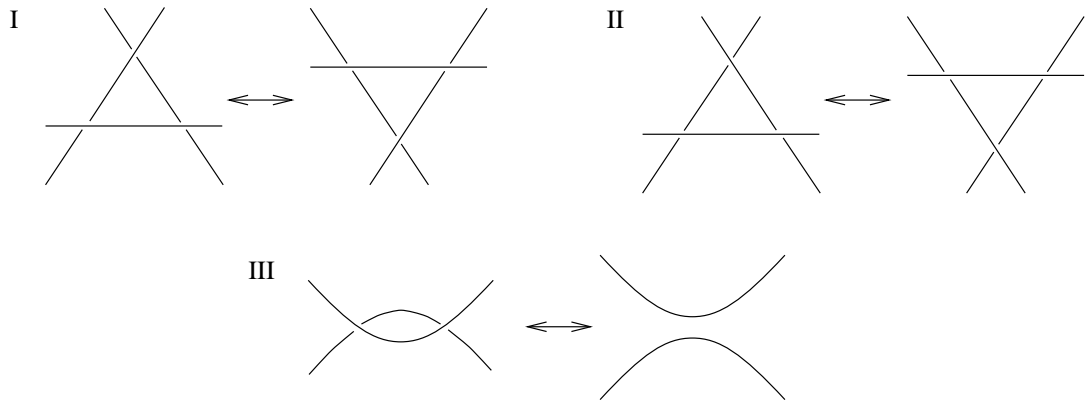}}
   \caption{The three Legendrian Reidemeister moves.}
   \mylabel{fig:reidemeister}
\end{figure}

\subsection{Move I}
\mylabel{ssec:moveI}

Let the algebras associated to the diagrams of $K$ before and after
move I be $\alg(a, b, c, v_1, \ldots, v_n; \df)$ and $\alg(a, b, c,
v_1,\ldots, v_n; \df')$ respectively, where the crossings $a$, $b$,
and $c$ are indicated in Figure~\ref{fig:moveI}.
In order to exhibit an elementary isomorphism between these two
DGAs, we need some more notation.  Figure~\ref{fig:moveI} labels
the twelve relevant quadrants by their orientation signs
$\varepsilon_{i,a},\varepsilon_{i,b},\varepsilon_{i,c}$.
Also, for either diagram in Figure~\ref{fig:moveI}, let
$\gamma_a,\gamma_b,\gamma_c$ be the capping paths corresponding to
the three crossings.  As we approach
the triple point intersection from either diagram,
$\gamma_a - \gamma_b - \gamma_c$ limits to a cycle in $H_1(K)$;
let $k$ be the number of times it winds around $K$.

\begin{figure}
  \centerline{\includegraphics{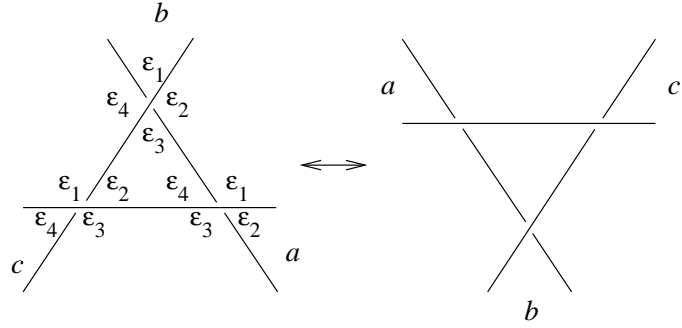}}
  \caption{A labeling for the orientation signs $\varepsilon_{i,a},
    \varepsilon_{i,b},\varepsilon_{i,c}$ for the twelve relevant
    quadrants in move I.  In the diagram, the subscript $a,b,c$ is
    suppressed.  }  \mylabel{fig:moveI}
\end{figure}

Define an elementary automorphism $\Phi$ on
$\alg=\alg(a,b,c,v_1,\ldots,v_n)$ by its action on the generators of
$\alg$:
\begin{equation} \mylabel{eqn:move1-isom}
  \Phi(w) = \begin{cases}
    a - \varepsilon t^{-k}cb & \text{if $w=a$} \\
    w & \text{otherwise,}
  \end{cases}
\end{equation}
where $\varepsilon =
\varepsilon_{4,a}\varepsilon_{4,b}\varepsilon_{1,c}$.  We claim that
$\Phi$ gives the desired tame isomorphism between the two DGAs.

We first note the following lemma, whose proof is similar to the proof
of Lemma~\ref{lem:nonconvexcorners} and is omitted.
\begin{lemma}
  $\sgn a = (\sgn b)(\sgn c)$.  \mylabel{lem:sumdegs}
\end{lemma}
Similarly, it is not hard to show that $\Phi$ preserves the grading on
$\alg$.

It remains to show that $\Phi$ is a chain map, intertwining the actions
of $\df$ and $\df'$.  We demonstrate this on the generators of $\alg$.
It suffices to consider $a$ and the generators whose differential has
a term containing $a$, since move I does not change any of the other
differentials (see \cite{chv}).

First consider a crossing $v$ for which $\df v$ contains $a$.
We may group the terms in $\df v$ and $\df' v$ as follows:
\begin{enumerate}
\item the terms in $\df v$ containing $cb$, and the corresponding terms
in $\df v$ and $\df' v$ replacing $cb$ by $a$;
\item
the terms in
$\df' v$ containing $cb$, and the corresponding terms in $\df v$
and $\df' v$ replacing $cb$ by $a$;
\item the remaining terms in
$\df v$ and $\df' v$, which are identical and do not contain $a$ or
$cb$.
\end{enumerate}
For $i=1,2,3$, denote by $\df_i v$ and $\df_i' v$ the
contributions of these three groups to $\df v$ and $\df' v$.

It is straightforward to check that $\df_1 v$ is simply $\df'_1 v$
with $a$ replaced by $a + \varepsilon_{4,a}\varepsilon_{4,b}\varepsilon_{1,c}t^{-k} cb$,
and so $\df'_1 v = \Phi \df_1 v$.  Similarly, $\df_2' v$ is
$\df_2 v$ with $a$ replaced by
$a + \varepsilon_{2,a}\varepsilon_{2,b}\varepsilon_{3,c} t^{-k} cb$.
By Lemmas~\ref{lem:prodsigns} and \ref{lem:sumdegs}, we have
\[
\varepsilon_{2,a}\varepsilon_{4,a}\varepsilon_{2,b}\varepsilon_{4,b}
\varepsilon_{1,c}\varepsilon_{3,c} = (-\sgn a)(-\sgn b)(-\sgn c) = -1;
\]
it follows that $\df'_2 v  = \Phi \df_2 v$.  Since we also
have the trivial identity $\df'_3 v = \Phi \df_3 v$, we conclude
that $\df' \Phi v = \df' v = \Phi \df v$.

Finally, we consider $\df a$ and $\df' a$.  We may group the terms
in $\df a$ and $\df' a$ as follows:
\begin{enumerate}
\item the disks in either $\df a$ or $\df' a$ with positive corner
in the quadrant labeled $\varepsilon_{1,a}$, which do not have an
adjacent corner at $b$ or $c$;
\item the disks in either $\df a$ or $\df' a$ with positive corner
at $\varepsilon_{3,a}$ and no adjacent corner at $b$ or $c$;
\item the disks in $\df a$ (resp.\ $\df' a$)
with positive corner at $\varepsilon_{1,a}$ (resp.\ $\varepsilon_{3,a}$)
and adjacent corner at $b$;
\item the disks in $\df a$ (resp.\ $\df' a$) with positive corner at
$\varepsilon_{3,a}$ (resp.\ $\varepsilon_{1,a}$)
and adjacent corner at $c$.
\end{enumerate}
As before, let $\df_i a$, $\df_i' a$ be the contributions of
these four groups to $\df a$, $\df' a$, for $i=1,2,3,4$.

Clearly $\df_1 a = \df_1'a$ and $\df_2 a = \df_2' a$.  On the other
hand, gluing the middle triangle to any disk in $\df_3 a$ or $\df_3'
a$ gives a disk with positive corner at $c$, and any disk in $\df c =
\df' c$ is obtained this way.  The sign and power-of-$t$ difference
between terms in $\df_3 a$ and the corresponding terms in $\df' c$ is
$\varepsilon_{1,a}\varepsilon_{2,b}\varepsilon_{2,c} t^{-k} =
-\varepsilon t^{-k}$ by Lemma~\ref{lem:prodsigns}, while the
difference between terms in $\df_3' a$ and the corresponding terms in
$\df' c$ is $\varepsilon_{3,a}\varepsilon_{4,b}\varepsilon_{4,c}
t^{-k} = \varepsilon t^{-k}$; hence $\df_3' a - \df_3 a = \varepsilon
t^{-k} (\df' c) b.$ Similarly, $\df_4' a - \df_4 a = (\sgn
c)\varepsilon t^{-k} c (\df' b)$.  We conclude that
\begin{equation}
  \df' a - \df a = \varepsilon t^{-k} ((\df' c) b + (\sgn c) c (\df' b))
  = \varepsilon t^{-k} \df'(cb).
\end{equation}
Hence $\Phi \df a = \df' \Phi a$, as desired.

\subsection{Move II}
\mylabel{ssec:moveII}

Let the algebras associated to the diagrams of $K$ before and after
move II be the same as for move I.  In this case, we claim that $\df =
\df'$.  The fact that there are no new disks, unlike in move I, is
derived from the following corollary of Stokes' theorem:

\begin{lemma}[Chekanov]
  \mylabel{lem:height} Let $u: (D^2_*, S^1_*) \to (\cc, K)$ be a
  holomorphic disk with positive punctures $t_1, \ldots, t_n$ and
  negative punctures $\tau_1, \ldots, \tau_m$.  Denote the height of
  the crossing $a$ by $h(a)$.  Then
  \begin{equation}
    \mylabel{eqn:heights}
    \sum_1^n h(u(t_i)) - \sum_1^m h(u(\tau_i)) = \int_{D^2_*} u^*(dx \wedge dy)
    > 0.
  \end{equation}
\end{lemma}

\noindent
See \cite{chv} for a detailed account.\footnote{Chekanov uses ``move
IIIa'' to denote our move II.}  The signs and powers of $t$
do not change since all
disks cover exactly the same quadrants at the crossings before and
after the move.  Hence $\df = \df'$.

\subsection{Move III}
\mylabel{ssec:moveIII}

Let $a$ and $b$ denote the two new crossings produced by move III,
as labeled in Figure~\ref{fig:moveIII}.  We may then write the
algebras associated to the diagrams of $K$ before and after
move III as
\begin{align*}
\alg &= \alg(a, b, a_1, \ldots, a_n, b_1, \ldots, b_m;\df) \\
\alg' & = \alg(a_1 \ldots, a_n, b_1, \ldots, b_m; \df').
\end{align*}
Furthermore, suppose that the other crossings are ordered by height:
\begin{equation}
  \mylabel{eqn:height-order}
  h(a_n) \geq \cdots \geq h(a_1) \geq h(a) > h(b) \geq h(b_1)
  \geq \cdots \geq h(b_m).
\end{equation}

The orientation signs
$\varepsilon_a,\varepsilon'_a,\varepsilon_b,\varepsilon'_b$ for four
relevant quadrants have been labeled in Figure~\ref{fig:moveIII}.
Note that, due to our choice of capping paths in
Section~\ref{ssec:alg}, $\gamma_a$ and $\gamma_b$ limit to the same
path as $a$ and $b$ approach each other; hence the term in $\df a$
corresponding to the $2$-gon with corners at $a$ and $b$ is
$\varepsilon_a \varepsilon_b b$.  Lemma~\ref{lem:height} then tells us
that $\df a = \varepsilon_a\varepsilon_b b + \varepsilon_a' v$, where
$v$ is a sum of terms in the $b_i$ and $t$.  Note for future reference
that Lemma~\ref{lem:prodsigns} implies that
$\varepsilon_a \varepsilon'_a = -\sgn a$ and $\varepsilon_b
\varepsilon'_b = -\sgn b$; also, since \df\ lowers degree by $1$, we have
$(\sgn a)(\sgn b) = -1$.

\begin{figure}
    \centerline{\includegraphics{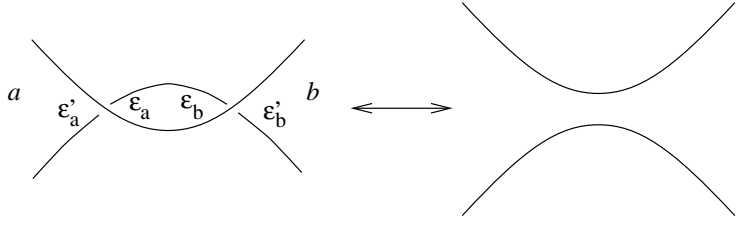}}
\vspace{12pt}
    \caption{A labeling for the crossings in move III;
      $\varepsilon_a,\varepsilon'_a,\varepsilon_b,\varepsilon'_b$ are
      orientation signs at crossings $a$ and $b$.}
    \mylabel{fig:moveIII}
\end{figure}

We now define a grading-preserving elementary isomorphism $\Phi_0:
\alg \to S_{|a|}(\alg')$ by its action on generators:
\begin{equation}
  \Phi_0(w) = \begin{cases}
    e_1 & w=a \\
    \varepsilon_a \varepsilon_b e_2 + \varepsilon'_b v & w=b \\
    w & \text{otherwise.}
  \end{cases}
\end{equation}
Although $\Phi_0$ is not a chain map on $\alg$, it is not far off, as
illustrated by the two lemmas below.  Define $\alg_i =
S_{|a|}(\alg(a_1, \ldots, a_i, b_1, \ldots, b_m))$, and let $\tau:
S(\alg) \rightarrow \alg$ be the obvious projection map.

\begin{lemma}
  \mylabel{lem:restr-chain} $\Phi_0|_{\alg_0}$ is a chain map.
\end{lemma}

\begin{proof}
  We prove this on generators of $\alg_0$.  There is nothing to prove
  for the $b_i$ since, by Lemma~\ref{lem:height}, $\df b_i$ contains
  only terms involving $b_j$ with $j>i$.  On the other hand, direct
  computation shows that $\Phi_0 \df a = e_2 = \df' \Phi_0 a$ and
  $\Phi_0 \df b = \df b = \df \Phi_0 b$.
\end{proof}

\begin{lemma}
  \mylabel{lem:proj-chain}
  $\tau \circ \df' \circ \Phi_0 = \tau \circ \Phi_0 \circ \df.$
\end{lemma}

\begin{proof}
  By Lemma~\ref{lem:restr-chain}, it suffices to prove equality on the
  generators $a_i$.  Let $W_1$ denote the sum of the terms which
  appear in both $\df a_i$ and $\df' a_i$; let $W_2$ denote the sum of
  the terms in $\df a_i$ involving $b$; and write $\df a_i = W_1 + W_2
  + W_3$ and $\df' a_i = W_1 + W_4$.

  The terms in $W_1$ do not contain $a$ or $b$, and so $\Phi_0 W_1 =
  W_1$.  The terms in $W_3$ must involve $a$; since $\Phi_0(a) = e_1$,
  we have $\tau \Phi_0 W_3 = 0$.

  Now consider the terms in $W_4$; these arise from disks of the type
  shown on the right in Figure~\ref{fig:moveIII-glue}.  There is a
  one-to-one correspondence between these disks and pairs of disks in
  $\df$, one with positive corner at $a_i$ and a negative corner at
  $b$, and one with positive corner at $a$.  Thus $W_4$ is the
  result of taking $W_2$ and replacing every occurrence of $b$ by
  $\varepsilon'_a \varepsilon'_b (\df a - \varepsilon_a \varepsilon_b
  b) = \varepsilon'_b v$; in other words, $W_4 = \tau \Phi_0 W_2$.

  We conclude that
  \begin{equation}
    \tau \df' \Phi_0 a_i = \df' a_i = W_1 + W_4
    = W_1 + \tau \Phi_0 W_2 = \tau \Phi_0 \df a_i,
  \end{equation}
  as desired.
\end{proof}

\begin{figure}
    \centerline{\includegraphics[angle=270]{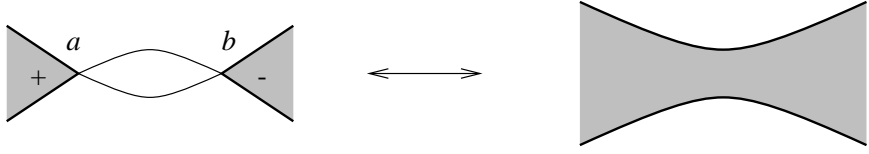}}
\vspace{12pt}
    \caption{Piecing together two disks from \df\ to get a disk in
    $\df'$.  The signs are Reeb signs, and the
crossing $a_i$ is schematically located off to the right
in both figures.}
    \mylabel{fig:moveIII-glue}
\end{figure}

Using Lemmas~\ref{lem:restr-chain} and \ref{lem:proj-chain}, we
bootstrap $\Phi_0$ up to a map $\Phi_n$ which is the desired chain map
on $\alg$, by inductively defining maps $\Phi_i$ which are chain maps
when restricted to $\alg_i$.  We define $\Phi_i$ along with elementary
automorphisms $g_i$ of $S(\alg)$ as follows.  Let $H:S(\alg') \to
S(\alg')$ be the map from the proof of Corollary~\ref{cor:homology}.
Define $g_i$ to fix all generators except $a_i$, and
\begin{equation} \mylabel{eqn:g}
  g_i(a_i) = a_i + H(\df' a_i - \Phi_{i-1} \df a_i);
\end{equation}
then define $\Phi_i = g_i \Phi_{i-1}$.

We collect several facts that we will need.  Recall
equation~(\ref{eqn:chain-homotopy}) from the proof of
Corollary~\ref{cor:homology}, namely:
$$\tau - Id = \df' H + H \df'.$$
Also, since $\tau H = 0$, we have
$\tau g_i = \tau$ for all $i$, and hence $\tau \Phi_i = \tau \Phi_0$
for all $i$.  Finally, note that, because the $a_j$ are ordered by
height, $\df a_i \in \alg_{i-1}$.

Now assume that $\Phi_{i-1}|_{\alg_{i-1}}$ is a chain map, i.e., that
$\df' \Phi_{i-1} = \Phi_{i-1} \df$ on $\alg_{i-1}$; we show that
$\Phi_i|_{\alg_i}$ is a chain map.  We can now calculate, at one point
using Lemma~\ref{lem:proj-chain}:
\begin{equation}
  \begin{split}
    \Phi_{i-1} \df a_i &= \tau \Phi_{i-1} \df a_i -
    \df' H \Phi_{i-1} \df a_i - H \df' \Phi_{i-1} \df a_i \\
    &= \tau \Phi_0 \df a_i - \df' H \Phi_{i-1} \df a_i
    - H \Phi_{i-1} \df^2 a_i \\
    &= \tau \df' \Phi_0 a_i - \df' H \Phi_{i-1} \df a_i \\
    &= (Id + \df' H + H \df') \df' a_i - \df' H \Phi_{i-1} \df a_i \\
    &= \df'(a_i + H\df' a_i - H \Phi_{i-1} \df a_i) \\
    &= \df' g_i(a_i).
  \end{split}
\end{equation}
Since $\Phi_{i-1} \df a_i \in \alg_{i-1}$, it follows that $\Phi_i \df
a_i = \df' g_i(a_i) = \df' \Phi_i a_i$.  On the other hand, the
induction hypothesis implies that $\Phi_i \df = \df' \Phi_i$ on
$\alg_{i-1}$.  Hence $\Phi_i$ is a chain map on $\alg_i$, completing
the induction.

This concludes the proof of Theorem~\ref{thm:invariance}. \qed

\bigskip
\begin{center}
{\bf
Part II \break
A Geometric Framework
}
\end{center}
\medskip

\noindent
In this part of the paper, we show that the combinatorial theory
developed in the previous sections fits into a much richer geometric
framework. In Section~\ref{sec:ch}, we show that our combinatorial
theory is a faithful translation of Eliashberg and Hofer's definition
of (relative) contact homology. Their theory is much more general and
provides invariants for Legendrian submanifolds in any contact
manifold. However, it is hard to make computations in their setup, so
it is helpful to know that the (easy) combinatorial definition of the
contact homology for Legendrian knots in $\rr^3$ is equivalent to
their definition. Moreover, having Eliashberg and Hofer's geometric
ideas in mind makes some of the combinatorial proofs more transparent.
Finally, in Section~\ref{sec:ori} we describe ``coherent
orientations'' in contact homology and show that, when translated into
our combinatorial framework, they yield the sign conventions described
in Part I.

\section{Relative Contact Homology}
\mylabel{sec:ch}

The goal of this section is to sketch a relative version of
Eliashberg and Hofer's contact homology theory.  Our presentation
is similar to that of \cite[Section 2.7]{egh}, in which relative
contact homology is set in the more general context of
``symplectic field theory''.  We will then specialize to the case
of Legendrian knots in the standard contact structure on $\rr^3$
and will show how to project the general theory into the $xy$
plane.  There the objects in the definition of contact homology
can be seen to be combinatorial in nature. Note that the analytic
details of the general theory have yet to be completed. When the
general theory is worked out, our translation shows that
computations in the combinatorial setting are also computations in
the general theory.

\subsection{The General Case}
\label{ssec:ch-general}

Let us begin by setting up some notation.  Let $(M; \xi)$ be a contact
3-manifold with contact form $\alpha$.  Let $K$ be a Legendrian knot
in $M$.  For simplicity of presentation, we consider only the case
where $H_1(M) = 0$ and $H_2(M) = 0$.  Thus, we may identify $H_2(M,
K)$ and $H_1(K)$ via the standard boundary homomorphism.  Further, a
choice of orientation on $K$ induces an identification of $H_1(K)$
with \zz. Let $X_\alpha$ be the Reeb vector field of $\alpha.$ A
segment of a flow line for $X_\alpha$ starting and ending on $K$ is
called a \textbf{Reeb chord} for $K.$

\begin{definition}
  \mylabel{defn:ch-alg} Let \alg\ be the free associative graded
  unital algebra over the group ring $\zz[H_2(M, K)] = \zz[H_1(K)]$
  generated by the Reeb chords in $(M,K; \alpha)$.  The generators of
  \alg\ are graded by their Conley-Zehnder indices (see below).  The
  generator $t$ of $\zz[H_1(K)]$ is graded by $2\,r(K)$.
\end{definition}

To define the Conley-Zehnder index of a Reeb chord $a(t)$, we must fix
a ``capping path'' $\gamma_a$ inside $K$ that connects $a(1)$ to
$a(0)$.  Next, we choose a surface $F_a$ with $\partial F_a = a \cup
\gamma_a$ and a trivialization of $\xi$ over $F_a$.  Let
$E$ be the sub-bundle of $\xi$ over $a \cup \gamma_a$ defined by:
\begin{equation}\mylabel{eqn:cz}
  \begin{split}
    E|_{\gamma_a(t)} & =  T_{\gamma_a(t)}K \\
    E|_{a(t)} & =  D\Phi_\alpha(t) \cdot T_{a(0)}K,
  \end{split}
\end{equation}
where $\Phi_\alpha$ is the flow of $X_\alpha.$
Using the trivialization, the sub-bundle $E$ may be viewed as a path
of Lagrangian subspaces in a fixed symplectic vector space.  Let
$CZ(a)$ be the Conley-Zehnder index of this path, as defined
in~\cite{robbin-salamon:maslov}.

While the Conley-Zehnder index is independent of the choices of $F_a$
and the trivialization of $\xi$ over $F_a$, it does depend on the
choice of capping path $\gamma_a$.  Suppose that $\tilde{\gamma_a}$ is
another such choice.  Since the paths $\gamma_a$ and
$\tilde{\gamma_a}$ have the same starting and ending points, they
differ up to homotopy by a path $\gamma_n$ that winds $n$ times around
the knot $K$.  We have:
\begin{equation} \mylabel{eqn:cz-ambiguity}
  \begin{split}
    CZ(\gamma) & =  \mu(\gamma_n) + CZ(\tilde{\gamma})  \\
    & =  n \, \mu (\gamma_1) + CZ(\tilde{\gamma}).
  \end{split}
\end{equation}
Here, $\mu$ is the Maslov index of a loop of Lagrangian subspaces.  As
we saw in Section~\ref{ssec:alg}, we can get a true \zz\ grading on
$A$ by making the choice of capping path explicit in the algebraic
structure.

To define a differential $\df$ on \alg, we must consider the
symplectization $(M\times\rr, \omega=d(e^\tau\alpha))$ of $(M, \xi)$,
and fix an almost complex structure $J$ which is compatible with the
symplectic form ({\em i.e.} $\omega(v, Jv)>0$ for $v\not=0$), and
which, in addition, satisfies:
\begin{equation} \mylabel{eqn:gen-J}
  \begin{split}
    J(\partial_\tau) & =  X_\alpha \\
    J(\xi) & =  \xi.
  \end{split}
\end{equation}
Now the differential \df\ on \alg\ is defined by counting certain rigid $J$-holomorphic
curves in the symplectization $(M \times \rr, d(e^\tau \alpha))$.

The $J$-holomorphic curves that we count are maps
\begin{equation} \mylabel{eqn:j-disks}
  f: (D^2_*, \partial D^2_*) \to (M \times
  \rr, K \times \rr)
\end{equation}
where $D^2_* = D^2 \setminus \{x, y_1, \ldots, y_n\}$ and $ \{x, y_1,
\ldots, y_n\}$ lies in $\partial D^2$.  Let $f_M$ be the projection of
$f$ to $M$ and let $f_\rr$ be the projection of $f$ to \rr.  We want
these disks to have boundary in the Lagrangian submanifold $K \times
\rr$ and to satisfy some asymptotic conditions near the punctures.
Note that a neighborhood of a
boundary puncture $x \in \partial D^2$ is conformally equivalent to the strip
$(0, \infty) \times [0,1]$ with coordinates $(s,t)$ such that
approaching $\infty$ in the strip is equivalent to approaching $x$ in
the disk.

\begin{definition}
  \mylabel{defn:reeb-converge} We say that $f$, parametrized as above
  near a boundary puncture $x$, \textbf{tends asymptotically to a Reeb
    strip} over the Reeb chord $a(t)$ at $\pm\infty,$ if:
  \begin{eqnarray*}
    \lim_{s \to \infty} f_{\rr}(s,t) & = & \pm \infty
    \mylabel{eqn:punc-limit1} \\
    \lim_{s \to \infty} f_{M}(s,t) & = & a(t) \mylabel{eqn:punc-limit2}.
  \end{eqnarray*}
\end{definition}

See~\cite{abbas} for convergence results for $J$-holomorphic curves
near boundary punctures.  We are now ready to define the moduli spaces
involved in the differential.

\begin{definition}
  \mylabel{defn:ms} $\ms^A(a;b_1,\ldots,b_n)$, the
  \textbf{moduli space of $J$-holomorphic disks}
  realizing the homology class $A \in H_2(M,K) \cong H_1(K)$, and
  with a positive puncture at the Reeb chord $a$ and (cyclically
  ordered) negative punctures at the Reeb chords $b_1, \ldots, b_n$,
  consists of all proper
  $J$-holomorphic maps $f$ as in equation (\ref{eqn:j-disks}), that
  satisfy the following conditions:
  \begin{enumerate}
  \item The map $f$ has finite energy:
    \begin{equation} \mylabel{eqn:energy}
      \int_{D^2_*} f^*d\alpha < \infty.
    \end{equation}
  \item The cycle $f_M(\partial D_*^2) \cup \gamma_a \cup -\gamma_{b_1}
    \cup \cdots \cup -\gamma_{b_n}$ represents the homology class $A$.
  \item Near $x$, $f$ tends asymptotically to a Reeb strip over $a$ at
    $+\infty$.
  \item Near $y_j$, $f$ tends asymptotically to a Reeb strip over
    $b_j$ at $-\infty$.
  \end{enumerate}
  Two maps $f$ and $g$ in $\ms^A(a;b_1, \ldots, b_n)$ are equivalent
  if there is a conformal map $\phi:D^2_* \to D^2_*$ such that $f = g
  \circ \phi$.
\end{definition}

Since $J$ is invariant under translation in the ``vertical'' $\tau$
direction, every $J$-holomorphic disk is part of a vertically
invariant family.  To alleviate confusion, we will denote the
dimension of $\ms^A$ by $k+1$, where the $1$ indicates dimension in
the vertically invariant $\tau$ direction.

The local structure of the moduli space $\ms^A(a;b_1, \ldots, b_n)$
near a map $f$ may be analyzed using the Implicit Function Theorem.
Consider the linearization of the \dbar\ operator at $f$:
$$D_f \dbar: \Omega^0(f^*T(M\times \rr)) \to \Omega^{0,1}(f^*T(M\times
\rr)).$$
Assuming sufficient genericity of the data and a suitable
Sobolev (or H\"older) setup (see~\cite{hzw:fredholm}, for example),
this map is surjective and Fredholm.  The Implicit Function Theorem
gives a local parametrization of $\ms^A(a;b_1, \ldots, b_n)$ whose
domain is a neighborhood of zero in the kernel of $D_f \dbar$.  The
dimension of the moduli space is given by:

\begin{proposition}[\cite{bourgeois:dim}]
  \mylabel{prop:dimension}
  \begin{equation*} \mylabel{eqn:gen-dim}
    \dim \ms^A(a;b_1, \ldots, b_n) = CZ(a) -\sum CZ(b_j) + 2r(K) \cdot A.
  \end{equation*}
\end{proposition}

Finally, we define a differential on \alg\ using all
$(0+1)$-dimensional moduli spaces:

\begin{definition}
  \mylabel{defn:ch-diff} Let $a, b_1, \ldots, b_n$ be Reeb chords.
  Define:
  \begin{equation*}\mylabel{eqn:d}
    \df(a) = \sum_{\dim (\ms^A(a; b_1, \ldots, b_n)) = 0+1} (\#\ms/\rr)
    t^A b_1 \cdots b_n.
  \end{equation*}
  Here, the number of points in \ms\ is counted with sign using
  an orientation on \ms\ described in Section~\ref{sec:ori}.
\end{definition}

In general, we expect that \df\ makes \alg\ into a differential graded
algebra whose homology is an invariant of the Legendrian knot and the
contact structure (see \cite{egh}). The remainder of this section
proves the invariance of $(\alg, \df)$ for Legendrian knots in $\rr^3$
with the fixed form $\alpha_0$ by translating everything into the
combinatorial setting of Part I.

\subsection{A Two-Dimensional Projection}
\label{ssec:proj}

We now specialize to the case of an oriented Legendrian knot $K$ in
$(\rr^3, \alpha_0)$. The first step in translating the general
relative contact homology theory into Chekanov's combinatorial DGA is
to project all of the objects involved in the definition of \alg\ into
the $xy$ plane.  Suppose that $K$ that admits a generic projection to
the $xy$ plane.

Since the Reeb field of $\alpha_0$ is $\partial_z$, the crossings of
the diagram of $K$ correspond to the Reeb chords of $(\rr^3, K,
\alpha_0)$.  We choose capping paths and grade the crossings by their
Conley-Zehnder indices.
Note that the contact planes $\xi_0$ may be globally trivialized by
$\partial_x$ and $\partial_y - x\,\partial_z$.  In the projection,
these are just the standard coordinate vector fields and hence, modulo
the correction at the end of the path, the Conley-Zehnder index is
just twice the rotation number of the path with respect to the
standard trivialization of $T\rr^2$.

Next, we describe an explicit $J$ for the symplectization $(\rr^3
\times \rr, d(e^\tau \alpha_0))$.  Recall that the compatible almost
complex structure described in Section~\ref{ssec:ch-general} must satisfy
equation (\ref{eqn:gen-J}).  Since the Reeb field is $\partial_z$, the
following $J$ works:
\begin{equation} \mylabel{eqn:r3-J}
  \begin{split}
    J(\partial_x) & =  \partial_y - x \partial_z \\
    J(\partial_y) & =  - x \partial_\tau - \partial_x \\
    J(\partial_z) & =  -\partial_\tau  \\
    J(\partial_\tau) & =  \partial_z.
  \end{split}
\end{equation}

Write a map $f:D^2_* \to \rr^3 \times \rr$ as
\begin{equation}
  f(u,v) = (x(u,v),y(u,v), z(u, v), \tau(u,v)).
\end{equation}
Using the $J$ in (\ref{eqn:r3-J}), we
can write down the $\dbar_J$ equations for $f$ as follows:
\begin{equation} \mylabel{eqn:dbar}
  \begin{split}
    \partial_u x - \partial_v y & =  0 \\
    \partial_u y + \partial_v x & =  0 \\
    \partial_u \tau - \partial_v z & =  x \partial_v y \\
    \partial_u z + \partial_v \tau & =  x \partial_u x.
  \end{split}
\end{equation}
Thus, the $xy$ projections of $J$-holomorphic maps of $D^2_*$ are
actually holomorphic as maps to \cc\ (i.e.\ the $xy$-plane endowed with
the standard complex structure).  It follows that every moduli space
of $J$-holomorphic disks projects to a family of holomorphic maps
\begin{equation} \mylabel{eqn:delta-maps}
  f:(D^2_*, \partial D^2_*) \to (\cc, \pi_{xy}(K)).
\end{equation}

The conditions on the $J$-holomorphic disks in $\ms^A(a; b_1, \ldots,
b_k)$ translate to the following restrictions on the maps $f$:
\begin{enumerate}
\item The homology class $\left[ \left( \pi_{xy}|_K \right) ^{-1}
    \left( \mathrm{Im} (f|_{\partial D^2}) \cup \gamma_a \cup
      -\gamma_{b_1} \cup \cdots \cup -\gamma_{b_n} \right) \right]$
  coincides with $A$.
\item The map $f$ sends the boundary punctures to the crossings of the
  diagram of $K$. At a boundary puncture, the map $f$ covers either
  one or three quadrants, with the majority of the quadrants positive
  at the crossing $a$ and negative at the crossings $b_i$.
\end{enumerate}

\begin{definition}
  \mylabel{defn:delta} The space of maps $\Delta_h^A(a;b_1, \ldots,
  b_n)$ consists of all holomorphic maps $f$ as in equation
  (\ref{eqn:delta-maps}) that satisfy conditions $1$ and $2$ above.

  We consider two maps $f$ and $g$ in $\Delta^A_h$ to be equivalent if
  there is a conformal map $\phi:D_*^2 \to D_*^2$ such that $f = g
  \circ \phi$.
\end{definition}
The spaces $\Delta^A_h(a;b_1, \ldots, b_n)$ are clearly subsets of the
combinatorially defined spaces $\Delta^A(a;b_1, \ldots, b_n)$
introduced in Definition~\ref{defn:delta-combinatorial}.  We use the
spaces $\Delta^A_h$ as a convenient intermediate step in showing the
equivalence between $\ms^A(a; b_1,\ldots, b_n)$ and
$\Delta^A(a;b_1,\ldots, b_n).$

Define a projection map $p$ by:
\begin{equation}
  \begin{split}
    \ms^A(a; b_1,\ldots, b_n) / \rr & \to  \Delta_h^A(a;b_1,
    \ldots, b_n) \\
    f & \mapsto \pi_{xy} \circ f.
  \end{split}
\end{equation}
Here the $\rr$-action is vertical translation.

\begin{theorem}[Translation Theorem]
  \mylabel{thm:translation} The following three spaces are homeomorphic:
  \begin{enumerate}
  \item $\ms^A(a; b_1,\ldots, b_n) / \rr$
  \item $\Delta_h^A(a;b_1, \ldots, b_n)$
  \item $\Delta^A(a;b_1, \ldots, b_n)$
  \end{enumerate}
  The projection $p$ gives an explicit homeomorphism between the first
  two.  The inclusion of $\Delta^A_h$ into $\Delta^A$ gives the
  homeomorphism between the last two.
\end{theorem}

Note that this implies:
\begin{equation} \mylabel{eqn:delta-ms-dim}
  \dim \Delta^A(a; b_1, \ldots, b_n) = \dim \ms^A(a; b_1, \ldots, b_n)
  -1.
\end{equation}

We will postpone the proof until the end of the section.  Since the
theorem tells us that $p$ is a homeomorphism, any system of
coordinates we find on the spaces $\Delta^A_h$ can be lifted to
coordinates on the moduli spaces $\ms^A$ (along with a coordinate that
parametrizes the vertical translations in $\ms^A$).  As it turns out,
we can use classical complex analysis to coordinatize $\Delta^A_h$.

\begin{proposition}
  \mylabel{prop:coords} The space $\Delta^A_h(a;b_1, \ldots, b_n)$ has
  local coordinates given by the images of $n$ interior branch points
  and $m$ branch points on the boundary. Thus,
  \begin{equation} \mylabel{eqn:delta-dim}
    2n+m=\dim \Delta^A_h.
  \end{equation}
\end{proposition}

\begin{proof}
  Let $f \in \Delta^A_h$.  Construct a Riemann surface $S$ for the
  inverse of $f$ by analytic continuation. Think of $S$ as a branched
  cover over the image of $f$.  Let $\tilde{f}$ be the lifting of $f$
  to a map from $D^2_*$ to $S$.  Since the inverse of $f$ is
  single-valued on $S$, $\tilde{f}$ must be an homeomorphism.  Using
  this construction, we must prove first that $f$ is the unique map
  in $\Delta^A_h$ with a given configuration of branch points in the
  image, and second that any small variation of the images of the
  branch points can be accomplished inside $\Delta^A_h$.

  For the first part, suppose that $g$ is an element of $\Delta^A_h$
  whose branch points have the same images as those of $f$.  Then $g$
  lifts as a biholomorphism to $\tilde{g}: D^2_* \to S$.  The map
  $\tilde{g} \circ \tilde{f}^{-1}$ is an automorphism of the disk, and
  hence $\tilde{g}$ and $\tilde{f}$ differ only by a
  reparametrization.  Projecting down to \cc, we see that $f$ and $g$
  differ only by a reparametrization, or in other words, $f \sim g$
  in $\Delta^A_h$.

  For the second part, let $S'$ be the Riemann surface obtained by
  perturbing the image of the branch locus of $S$ in \cc. Note that
  $S'$ projects to a different region in \cc \,when a boundary branch
  point is perturbed along the diagram of $K$.  Since small
  perturbations of the image of the branch locus do not change the
  topology of $S$, the Uniformization Theorem provides a
  biholomorphism $g: D^2 \to S'$ that projects to an element of
  $\Delta^A_h$.

  Finally, each interior branch point has two degrees of freedom,
  whereas a boundary branch point contributes but one; the formula
  (\ref{eqn:delta-dim}) follows.
\end{proof}

We end this section with a proof of Theorem~\ref{thm:translation}.

\begin{proof}[Proof of Theorem~\ref{thm:translation}]
  We prove this theorem in two steps. We show first that
  $\Delta^A_h=\Delta^A$, and then that the projection from $\ms^A /
  \rr$ to $\Delta^A_h$ is a homeomorphism.

  To prove $\Delta^A_h=\Delta^A$, we clearly need only to show that
  $\Delta^A\subset\Delta^A_h.$ To this end, let $f:D^2_* \to \cc$ be
  an immersion satisfying the conditions for $\Delta^A$.  Use $f$ to
  pull back the complex structure from $\cc$. Then the Riemann Mapping
  Theorem provides a conformal equivalence $g$ between the new complex
  structure and the standard structure on the interior of $D^2$.
  Hence $f \circ g$ is holomorphic on the interior of $D^2_*$.  That
  $f$ lies in $\Delta^A_h$ now follows from the proof of
  Proposition~\ref{prop:coords} and the fact that holomorphic maps
  preserve orientation.

  We are left to show that $p:\ms^A / \rr\to \Delta^A_h$ is a homeomorphism.
  It is clear from the discussion above that the image of $p$ lies in
  $\Delta_h^A$.  The interesting part of the proof lies in the
  construction of an inverse $q$ to $p$ that lifts maps in $\Delta_h^A$
  to $\ms^A$.  To do this, note that a few simple manipulations of the
  $\dbar_J$ equations (\ref{eqn:dbar}) give, for $g=(x,y,z,\tau) \in
  \ms^A$:

  \begin{lemma}
    \mylabel{lem:harmonic} $z(u,v)$ is harmonic.
  \end{lemma}

  Thus, given $f = (x,y) \in \Delta_h^A$, we define $z(u,v)$ by solving
  the Dirichlet problem with a (discontinuous) boundary condition
  which may be formulated as follows: suppose that $(u,v) \in \partial
  D_*^2$.  By assumption, $(x(u,v), y(u,v))$ lies in the diagram of $K$.
  Away from the crossings, let $z(u,v)$ be the $z$ coordinate of the
  knot $K$ that lies above $(x(u,v), y(u,v))$.  This defines the
  boundary condition $z(u,v)$ uniquely on $\partial D^2_*$.

  Once we have determined $z(u,v)$, a similar manipulation of the
  $\dbar_J$ equations (\ref{eqn:dbar}) yields:

  \begin{lemma} \mylabel{lem:curl}
    $\partial_{uv}\tau = \partial_{vu}\tau.$
  \end{lemma}

  Combined with the Poincar\'e Lemma, Lemma~\ref{lem:curl} tells us
  that we can find a $\tau(u,v)$ (unique up to an additive constant)
  so that $q(f) = (x,y,z,\tau)$ lifts $f$ and solves the $\dbar_J$
  equations on the interior of $D^2$.  In fact, the solutions extend
  continuously to the boundary away from the punctures $x, y_1,
  \ldots, y_n$.

  The lift $q(f)$ clearly satisfies condition $2$ of the definition of
  $\ms^A(a;b_1, \ldots, b_n)$.  The first and last two conditions,
  namely that of finite energy and of asymptotic approach to Reeb strips
  at the punctures, need proof.  We tackle the last condition first.

  We begin by describing a local model for the lifting $q(f)$ near a
  positive puncture in the special case where the crossing that is the
  image of $f(x)$ is bounded by the $x$ and $y$ axes.  Further,
  suppose that the lift $\tilde{L}_0$ of the $x$ axis has constant $z$
  coordinate $0$ and that the lift $\tilde{L}_1$ of the $y$ axis has
  constant $z$ coordinate $\pi/2$.  In the $xy$ plane, the exponential
  map takes the strip $\Sigma=\rr_+ \times i[0,\pi/2]$ to the lower
  right-hand quadrant. Consider the following lifting $q(f)$ of $f$:
  \begin{equation}\mylabel{eqn:straight-line}
    \begin{split}
      \Sigma & \to  \rr^3 \times \rr \\
      u+iv & \mapsto  \left(e^{-u}\cos(v), -e^{-u}\sin(v), v, u +
        \frac{e^{-2u}}{2}\cos^2 (v)\right).
    \end{split}
  \end{equation}
  It is straightforward to check that this map is $J$-holomorphic and
  tends asymptotically to the Reeb chord over the origin as $u$ goes to
  $\infty$.

  It is not hard to generalize this model to the case where
  $\tilde{L}_0$ and $\tilde{L}_1$ are arbitrary straight lines whose
  projections $L_0$ and $L_1$ pass through the origin.  Note that, in
  this case, the original map $f$ is $f(u+iv) = c_0 e^{-\nu(u+iv)}$,
  where $c_0 \in L_0$, $c_0 e^{i\frac{\pi}{2}\nu_0} \in L_1$, $0 < \nu_0 <
  1$, and $\nu = k - \nu_0$ for some integer $k$.

  A theorem of Robbin and Salamon shows that every positive puncture
  is $O(e^{-Ks})$-close to a straight-line model in the $xy$-plane:

  \begin{theorem}[Robbin and Salamon \cite{robbin-salamon:hol}]
    \mylabel{thm:rs}  Let $L_0$ and $L_1$ be curves in \cc\ that pass
    through the origin.  Let $f: \Sigma \to \cc$ be
    a holomorphic map that satisfies:
    \begin{enumerate}
    \item $f(\rr_+ \times \{\frac{\pi j}{2}\}) \subset L_j$ for $j=0,1$.
    \item $\lim_{u \to \infty} f(u,v) = \lim_{u \to \infty}
      \partial_u f(u,v) = 0$ uniformly in $v$.
    \end{enumerate}
    Then there exist constants $c_0 \in \cc^*$, $\nu \in \rr_+$, and
    $\delta \in \rr_+$ such that
    \begin{equation*}
     f(u,v) = c_0 e^{-\nu(u+iv)} + O(e^{-(\nu + \delta)u}),
    \end{equation*}
    with $c_0 \in T_0 L_0$, $c_0 e^{i\frac{\pi}{2}\nu_0} \in T_0 L_1$, $0 < \nu_0 <
    1$, and $\nu = k - \nu_0$ for some integer $k$.
  \end{theorem}

  Given a crossing with curves $\tilde{L}_0$ and $\tilde{L}_1$, let
  $f:\Sigma \to \cc$ be a holomorphic curve that satisfies the
  hypotheses of Theorem~\ref{thm:rs}.  Let $f_0$ be a straight-line
  solution with respect to the lines $T_0 \tilde{L}_0$ and $T_0
  \tilde{L}_1$.  The theorem asserts that these two solutions differ
  by $O(e^{-(\nu + \delta)u})$.  Further, a simple calculation using
  the fact that the $z$ coordinates of $\tilde{L}_j$ are determined
  from their $xy$ coordinates shows that the boundary conditions for the
  Dirichlet problem differ by $O(e^{-(\nu + \delta)u})$.  By the
  maximum principle, the $z$ liftings of $f$ and $f_0$ differ by at
  most $O(e^{-(\nu + \delta)s})$.  Furthermore, by explicitly writing
  out the formula for the liftings to the $\tau$ coordinate that come
  from the proof of the Poincar\'e Lemma (see \cite{bott-tu}, for
  example), we see that the difference there is $O(ue^{-(\nu +
    \delta)u})$.

  Thus, on the strip $\Sigma$, the lifting of the straight-line model
  is exponentially close to the lifting of the general case.  Since
  the straight-line case tends asymptotically to a Reeb chord, so must
  the general case.

  Finally, we have to show that $q(f)$ has finite energy.  Let $S_0$
  be a small half-disk around the puncture $x$ in the disk
  $D^2_*$.  Similarly, for each $i$, let $S_i$ be a half-disk
  around $y_i$.  Let $\Gamma$ be the boundary of $D^2 \setminus \cup
  S_i$. Stokes' Theorem gives:
  \begin{equation} \mylabel{eqn:stokes}
    \int_{D^2 \setminus \cup S_i} q(f)^*d\alpha  = \int_\Gamma
    q(f)^*\alpha.
  \end{equation}
  As the $S_i$ get smaller, $\Gamma$ approaches $\partial D^2_*$ and
  so, using the straight-line model and Theorem~\ref{thm:rs}, we have:
  \begin{equation} \mylabel{eqn:energy-converge}
    \int_\Gamma q(f)^*\alpha \to \int_a \alpha - \sum_i \int_{b_i}
    \alpha < \infty.
  \end{equation}

  This completes the proof that the inverse map $q: \Delta^A_h \to \ms^A
  / \rr$ is well-defined.  Using the same methods as for the
  characterization of the asymptotic behavior of $q(f)$, it is clear
  that $q$ is a continuous inverse to $p$. The theorem follows.
\end{proof}

\begin{remark}
  Theorem~\ref{thm:translation} and Proposition~\ref{prop:coords}
  apply to the more general setting in which disks with multiple positive
  boundary punctures are present. Such disks must be understood when
  trying to generalize contact homology to symplectic field theory.
\end{remark}

\section{Coherent Orientations}
\mylabel{sec:ori} In this section we describe Floer and Hofer's idea
of coherent orientations \cite{floer-hofer} in the context of contact
homology \cite{egh}.
This allows the signed count of points in $\ms^A$ used in the
definition of $\df.$ In Section~\ref{ssec:comb-ori}, we will then
translate these orientations on $\ms^A$ to orientations on $\Delta^A$
using Theorem~\ref{thm:translation}. We show in
Section~\ref{ssec:algorithm} that these orientations yield the sign
conventions used in Part I.

\subsection{Geometric Ideas}
\mylabel{ssec:geom-ori}

The geometric idea behind orienting the moduli spaces $\ms^A$
discussed in Section~\ref{sec:ch} comes from Floer and Hofer's
\textbf{coherent orientations}.  In \cite{floer-hofer}, they detailed
a program for orienting the moduli spaces relevant to Floer homology
for periodic orbits of Hamiltonian systems, in which the moduli spaces
involve maps of infinite cylinders that limit to periodic orbits at
each end.  In \cite{egh}, Eliashberg, Givental, and Hofer generalized
and refined the coherent orientation idea to the setting of symplectic
field theory, of which contact homology is a special case.  In this
section, we will adapt their ideas to the relative setting in $\rr^3$.

This first step in the coherent orientation program is to consider the
operators involved in the definition of the moduli space $\ms^A(a;b_1,
\ldots, b_k)$ rather than just the moduli space itself.  Recall that
$\ms^A(a;b_1, \ldots, b_k)$ is the space of $J$-holomorphic maps from
the boundary-punctured disk $D^2_* = D^2 \setminus \{x, y_1, \ldots,
y_n\}$ to $\rr^3 \times \rr$ that send the boundary of the punctured
disk to $K \times \rr$ and tend asymptotically to the Reeb chords $a,
b_1, \ldots, b_n$ at the punctures.  Let $f \in \ms^A(a;b_1, \ldots,
b_k)$.
The linearized operator
$D_f\dbar_J$ on $f^*T(\rr^3 \times \rr)$ is Fredholm in the proper
analytic setup.  When $D_f\dbar_J$ is surjective, the Implicit
Function Theorem gives a local coordinate system on $\ms^A(a;b_1,
\ldots, b_k)$.  Thus, in order to orient $\ms^A(a;b_1, \ldots, b_k)$,
it suffices to orient the kernel bundle of the operators $D_f\dbar_J$
for $f \in \ms^A$.

Instead of doing this, Floer and Hofer's idea was to expand the set of
operators under consideration beyond the $D_f\dbar_J$ operators to the
space of \textbf{Cauchy-Riemann-type} operators that share the
boundary conditions and asymptotic behavior with the ``honest''
operators that come from the moduli space.  In fact, their program
orients Cauchy-Riemann-type operators on all complex bundles over all
closed Riemann surfaces.

In order to define Cauchy-Riemann-type operators in the relative case,
we introduce holomorphic coordinates $(0, \infty) \times
i[0,\pi/2]$ near each puncture on the boundary of a Riemann surface
$S_*$ with boundary.  Compactify each puncture $y_k$ with an interval
$I_k$ to get $(0, \infty] \times i[0,\pi/2]$.  Call the compactified
Riemann surface $\hat{S}$.

\begin{definition}[Compare with \cite{egh}]
  \mylabel{defn:bundle}A \textbf{smooth complex vector bundle $E$ over
    $S_*$} is a smooth complex bundle over $\hat{S}$ together with
  Hermitian trivializations $\Phi_k: E|_{I_k} \to [0,\pi/2] \times
  \cc^n$ and a real sub-bundle $E^\partial$ over $\partial S_*$.  An
  isomorphism between bundles must respect the Hermitian
  trivializations over $I_k$ and the real sub-bundle $E^\partial$.
\end{definition}

In this setup, the appropriate sections of $E$ to consider are those
which map the boundary of $S_*$ to the real sub-bundle $E^\partial$.
In fact, we want to consider sections in the Sobolev space
$H^{1,2}_{loc}(S_*, E)$ that, in addition, lie in $H^{1,2}\left( (0,
  \infty] \times i[0,\pi/2], E \right)$ with respect to local
``strip'' coordinates near each puncture.  We call the space of such
sections $H^{1,2}(E)$; a similar definition applies to $L^2$ sections.

Given a bundle $E$, define $X_E$ to be the bundle whose fiber over $z$
consists of $(i, J)$-anti-linear maps $\varphi:T_z \hat{S} \to E_z$.
We want to consider operators $L:H^{1,2}(E) \to L^2(X_E)$ which have
the form:
\begin{equation} \mylabel{eqn:c-r-operator}
  (Lh) \cdot X = \nabla_X h + J\nabla_{iX}h + A(h) \cdot X.
\end{equation}
Here, $h$ is a section in $H^{1,2}(E)$, $X$ is a vector field on
$S$, and $A$ is a section of $\mathrm{Hom}_{\rr}(E, X_E)$.  The
operators $L$ must satisfy an asymptotic condition near the punctures.
Namely, in local strip coordinates $(s,t) \in (0, \infty] \times
i[0,\pi/2]$, $L$ must take the form:
\begin{equation} \mylabel{eqn:c-r-infinity}
  (Lh) \cdot \partial_s = \partial_s h - A(s,t)\cdot h.
\end{equation}
Furthermore, $A(s,t) \to -i\partial_t - a(t)$, where $a(t)$ is a smooth
path of self-adjoint operators on $L^2([0, \pi/2], \cc^n)$.  This
completes the definition of a Cauchy-Riemann-type operator in the
relative case.

In the correct functional analytic setup, all Cauchy-Riemann-type
operators are Fredholm and hence the space $\Sigma(E)$ of
Cauchy-Riemann-type operators on a given bundle $E$ with fixed
asymptotic data has a well-defined determinant bundle. (See
\cite{floer-hofer} for more details about the definition of a
determinant bundle.) Assume for now that all of the spaces
$\Sigma(E)$ are orientable, i.e.\ that their determinant bundles are
trivializable.  Denote an orientation on $\Sigma(E)$ by $\sigma(E)$.
Since we assume that $\Sigma(E)$ is orientable, it is enough to orient
the determinant line over one operator $L \in \Sigma(E)$ to determine
$\sigma(E)$.  As a result, we will frequently abuse notation and refer
to $\sigma(L)$ instead of $\sigma(E)$.

Following \cite{egh}, we want to put a \textbf{coherent set of
  orientations} on all spaces of Cauchy-Riemann-type operators on
bundles over Riemann surfaces with boundary.  Such a set of
orientations must satisfy three axioms:

\begin{description}
\item[Disjoint Union] Given bundles $E_j \to S_j$ for $j=1,2$,
  define the \textbf{disjoint union bundle} $E_1 \sqcup E_2 \to S_1
  \sqcup S_2$ by $(E_1 \sqcup E_2)|_{S_j} = E_j$.  If we have
  operators $L_j \in \Sigma(E_j)$, then there is an operator $L$ on
  $E_1 \sqcup E_2$ such that $L|_{H^{1,2}(E_j)} = L_j$.  The
  determinant line $\det L$ is canonically isomorphic to $\det L_1
  \otimes \det L_2$, and hence the orientations $\sigma(L_j)$ induce
  an orientation $\sigma(L_1) \otimes \sigma(L_2)$ on $\det L$.  The
  disjoint union axiom tells us to use this orientation:
  \begin{equation} \mylabel{eqn:disj-union}
    \sigma(L_1) \otimes \sigma(L_2) = \sigma(L).
  \end{equation}

\item[Direct Sum] Given bundles $E,F \to S$ and operators $L \in
  \Sigma(E)$ and $K \in \Sigma(F)$, there is a canonically defined
  operator $L \oplus K \in \Sigma(E \oplus F)$.  Once again,
  $\sigma(L)$ and $\sigma(K)$ induce an orientation $\sigma(L) \oplus
  \sigma(K)$ on $\det (L \oplus K)$.  The direct sum axiom states
  that:
  \begin{equation} \mylabel{eqn:direct-sum}
    \sigma(L) \oplus \sigma(K) = \sigma(L \oplus K).
  \end{equation}

\item[Cut and Paste] Let $S$ be a disjoint union of one or more
  Riemann surfaces with punctures on the boundary.  Let $E \to S$ be a
  bundle.  Let $\gamma_j: [0,1] \to S, j=1,2$ be real-analytic
  embeddings with disjoint images; if $\gamma_j(0) \neq \gamma_j(1)$,
  then both $\gamma_1$ and $\gamma_2$ must map $\{0,1\}$ to $\partial
  S$.  Let $\Phi: E|_{\gamma_1} \to E|_{\gamma_2}$ be a vector bundle
  isomorphism covering $\gamma_2 \circ \gamma_1^{-1}$.

  If we cut $S$ along the curves $\gamma_j$, we obtain a surface
  $\bar{S}$ with corners that has (possibly additional) boundary
  components $\gamma_j^\pm$.  The bundle $E$ gives rise to a bundle
  $\bar{E} \to \bar{S}$.  We can identify sections of $E$ with
  sections $\bar{h}$ of $\bar{E}$ that satisfy the boundary condition
  \begin{equation}\mylabel{eqn:c-p-boundary}
    \bar{h}|_{\gamma_j^-} = \bar{h}|_{\gamma_j^+} \quad \mathrm{for\ }
    j=1,2.
  \end{equation}
  In addition, the operator $L$ induces an operator $\bar{L}$ that
  acts on the sections of $\bar{E}$ that satisfy the boundary
  conditions (\ref{eqn:c-p-boundary}).

  If we shuffle the boundary conditions so that they read:
  \begin{equation} \mylabel{eqn:c-p-boundary-2}
    \begin{split}
      \Phi \bar{h}|_{\gamma_1^-} & =  \bar{h}|_{\gamma_2^+} \\
      \bar{h}|_{\gamma_2^-} & =  \Phi \bar{h}|_{\gamma_1^+},
    \end{split}
  \end{equation}
  then the sections satisfying (\ref{eqn:c-p-boundary-2}) correspond
  to sections of a new bundle $E'$ which is formed by identifying
  $\gamma_1^+$ to $\gamma_2^-$ and $\gamma_1^-$ to $\gamma_2^+$ on the
  base level and using $\Phi$ to identify the bundles at the fiber
  level.  The operator $L$ also induces a new operator $L'$ on the
  sections of $E'$.

  The boundary conditions (\ref{eqn:c-p-boundary}) and
  (\ref{eqn:c-p-boundary-2}) may be connected by a path in the space
  of bundles with mixed Riemann- and Riemann-Hilbert-type boundary
  conditions.  These boundary conditions induce a path of
  Cauchy-Riemann-type operators, whose determinant lines may be
  oriented by continuation.  Thus, via this path, the operator $L'$
  gets an orientation $\sigma(L, \gamma_j, \Phi)$ from $L$.  The cut
  and paste axiom states that:
  \begin{equation}
    \sigma(L, \gamma_j, \Phi) = \sigma(L').
  \end{equation}
\end{description}

See \cite{egh} for more details, especially with regards to the cut
and paste operation.  The cut and paste operation in \cite{egh} was
designed to replace the more limited gluing construction of
\cite{floer-hofer}.  What is important for us is that gluing is
retained as a special case of cutting and pasting. To see this,
suppose that operators in $\Sigma(E\to S_1)$ tend asymptotically to
$A_\infty$ at $+\infty$ near a boundary puncture $x$ and that
operators in $\Sigma(F \to S_2)$ tend asymptotically to $A_\infty$ at
$-\infty$ near a boundary puncture $y$.  By a homotopy, we may assume
that there is an operator $L \in \Sigma(E)$ that is constant for all
$s > \rho -1$ in local strip coordinates near $x$; similarly, assume
that there is an operator $K \in \Sigma(F)$ that is constant for all
$s < - \rho +1$ in local strip coordinates near $y$.  Choose paths
$\gamma_1(t) = (\rho -1, t)$ near $x$ and $\gamma_2(t) = (-\rho+1, t)$
near $y$.  The cut and paste operation results in an operator $L \# K$
on the ``boundary connected sum'' of $S_1$ and $S_2$ that is equal to $L$
on the old $S_1$ and $K$ on the old $S_2$ as well as a constant
operator on a strip.  The operator $L \# K$ is precisely the one
obtained from gluing in \cite{floer-hofer}.

In \cite{egh}, Eliashberg, Givental, and Hofer prove that there exists
a coherent set of orientations on \emph{closed} Riemann surfaces (with punctures).
Further, this set of orientations is determined by the three axioms
above and a choice of orientations on $\Sigma(S^2 \times
\cc)$, $\Sigma(\mathcal{O}(1) \to S^2)$, and the operators associated to
the $+\infty$ asymptotic data at punctures and the trivial vertical cylinder
(or all asymptotic data at punctures);
see \cite{egh}.  For the relative case to work, we must also
include an orientation on
$\Sigma(E=D^2 \times \cc, E^\partial=\partial D\times \rr)$
and the operators associated to asymptotic data at positive boundary punctures (as in \ref{eqn:c-r-infinity})
and the trivial vertical strip.
However, the relative case does not work in general, since the spaces
of Cauchy-Riemann-type operators --- in particular, their boundary
conditions --- are not necessarily orientable.

We now specialize to the case of orientations for the relative contact
homology of Legendrian knots in $(\rr^3, \xi_0)$.  Since Reeb field
for the standard contact structure is invariant under translations, it
should come as no surprise that every Reeb chord yields the same
asymptotic data.  More precisely, we get:

\begin{lemma}
  \mylabel{lem:asymptotic-data} All asymptotic operators for $(\rr^3,
  \xi_0)$ have the form
  $$A_\infty h(t) = -i \frac{\partial h}{\partial t}$$
  with respect to
  the Hermitian trivialization $\{\partial_x, \partial_y +
  x\partial_z, \partial_\tau, \partial_z\}$.
\end{lemma}

This lemma follows from the definitions, the local model discussed in
Section~\ref{ssec:proj}, and Theorem~\ref{thm:rs}.

With this in mind, we set up some notation: the space of all
Cauchy-Riemann-type operators on a $\cc^2$-bundle over the
boundary-punctured disk $D^2_*$ with asymptotic data given by the Reeb
chords $\{a, b_1, \ldots, b_n\}$ is denoted $\Sigma(a; b_1, \ldots,
b_n)$.  When $\Sigma$ comes from an honest moduli space, the bundle
$f^*(T\rr^3 \times \rr)$ splits into $\cc \oplus \cc$, with the first
factor coming from the contact structure and the second coming from
the complex direction spanned by the Reeb orbit and the vertical
direction.  The boundary conditions for the first factor come from the
diagram of the knot $K$; we may regard them as a path $\Lambda$ in
$\rr \pp^1$ with fixed endpoints.  The boundary conditions are trivial
for the second factor as they always point in the vertical direction.
If $\Sigma(a; b_1, \ldots, b_n)$ has such a splitting, then it is
parametrized by triples $(J, A, \Lambda)$.  Each component of this
space of triples is contractible,\footnote{However, in
  higher dimensions, $\pi_2(\mathrm{Lagr}(\cc^n))$ is non-trivial, and
  hence the space of triples is not contractible.}  so the determinant
bundle over $\Sigma(a; b_1, \ldots, b_n)$ is orientable.

In the next section, we will describe a system of coherent
orientations for the moduli spaces involved in relative
contact homology for the standard contact $\rr^3$.  We will not
present a complete system; rather, we will only use the axioms to
orient the spaces required for the definition of the differential \df\
on \alg.  In particular, we only use split bundles as above, so the
spaces $\Sigma$ will always be orientable.  Using the Translation
Theorem (\ref{thm:translation}), we will translate the geometric
considerations detailed above into a diagrammatic formulation of
coherent orientations in the $xy$-plane.

\subsection{A Combinatorial Approach}
\mylabel{ssec:comb-ori}

Recall that the Translation Theorem~\ref{thm:translation} asserts that all
$J$-holomorphic curves in the symplectization $\rr^3 \times \rr$ can
be represented as holomorphic maps of boundary-punctured disks into
$\rr^2$ whose boundary is sent into the diagram of $K$. Slightly
rephrased, this means that the boundary of a $J$-holomorphic curve in
the symplectization gives rise to an oriented loop in the diagram of $K$
that is immersed except at finitely many
points, which may be either branch points or corners.
Any one loop encodes all of the boundary conditions (via pulling back
$TK$) and asymptotic conditions (via local lifts near the corners) for
a Cauchy-Riemann-type operator on a $\cc^2$-bundle over $D^2_*$; in
other words, a family of oriented loops is equivalent to a family of operators
$\Sigma(a; b_1, \ldots, b_n)$. Following Floer and Hofer, we will
expand the families of loops under consideration and orient \emph{all}
families $\Sigma(a; b_1, \ldots, b_n)$ of oriented loops $\gamma(s)$ in the
diagram of $K$ that are immersed except at finitely many
points, which may be either branch points or corners.

We put coordinates on a family of loops with $k$ branch points by
numbering the branch points consecutively with respect to the given
orientation of the loop.  We view them as coordinates on the $k$-torus
$\prod_k S^1$ (in general, these coordinates only lie in a subspace of
the $k$-torus). For orientation purposes, two numberings
are equivalent if they differ by an \emph{even} cyclic permutation.
An orientation on a loop with $k$ boundary branch points is given by
an ordered basis for $T_p (S^1)^k \oplus \rr$, which we shall
represent by the ordered set of vectors $\left< \pm v_1, \ldots, \pm v_k,
  \pm \partial_\tau \right>$ with $v_j \in T_p S^1$ agreeing with the orientation
on the loop $\gamma$ and $\partial_\tau \in
T_p \rr$.

Note that we may arbitrarily add a pair of branch points on a loop in
$\Sigma(a; b_1, \ldots, b_n)$ and still have a loop in $\Sigma(a; b_1,
\ldots, b_n).$ We will see below that the orientation is unaffected by
such an addition. This is related to the fact that if a loop $\gamma$
comes as the boundary of a $J$-holomorphic disk with $n$ interior
branch points and $m$ boundary branch points then, as we move about
the moduli space containing this disk, one of the interior branch
points could migrate to the boundary where it spawns two boundary
branch points. However, since the interior branch points have two
degrees of freedom with a canonical complex orientation, we have a
natural orientation associated to the two new boundary branch points.


\begin{figure}
  \centerline{\includegraphics{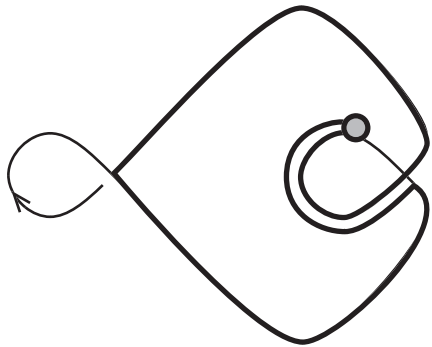}}
  \caption{The $xy$ projection of a $(1+1)$-dimensional moduli space.
    In this figure, and in all subsequent figures, the diagram of the
    knot $K$ is designated by a thin line, while the image of an
    element in a moduli space is designated by a thick line.}
  \mylabel{fig:xy-ms}
\end{figure}

The construction of a coherent system of orientations proceeds in
three steps:
\begin{enumerate}
\item Orient all families of loops with no corners.
\item At each crossing, orient precisely one family of loops with a
  single corner and no branch points with $\left<\partial_\tau
  \right>$.
\item Every other family of loops gets oriented by gluing to families
  of loops from step $2$ until there are no remaining corners and then
  comparing the result to the orientations from step $1$.
\end{enumerate}

To apply the cut-and-paste axiom we need to first choose an
orientation for the trivial vertical strip. This is a dimension $0$
object and will be oriented by $+1.$ Now using the cut-and-paste
axiom, the first step is equivalent to choosing orientations for
operators on the trivial bundle $D^2 \times \cc$ and on
$\mathcal{O}(1)$.  The second step mirrors a choice of orientation on
a single asymptotic operator.  See \cite{egh} for a similar procedure
in the closed case.

\textbf{Step 1.}  Any loop with no branch points or corners is related
to an operator on a bundle that can be obtained from $D^2\times \cc$
and $\mathcal{O}(1)$ using the cut-and-paste axiom. We choose these
orientations so that such a loop is oriented by
$\left<\partial_\tau\right>.$ Examining our discussion above
concerning the migration of interior branch points to the boundary
yields the following convention for a loop $\Sigma$ with $k$ branch
points and no corners:
label each segment between two branch points in $\Sigma$ by a
`$+$' if the orientations of $\Sigma$ and the knot $K$ agree and
by a `$-$' otherwise.  We will refer to the order of signs before
and after a branch point as its \textbf{alignment}.  Starting with
a branch point whose alignment is a $+$ then a $-$, number each
branch point consecutively in the order encountered; this gives an
ordering of the coordinates on $\Sigma$.  The orientation is given
by:
\begin{equation} \mylabel{eqn:no-corner-ori}
  \left<v_1, v_2, \ldots, v_{k-1}, v_k, \partial_\tau \right>.
\end{equation}
This completes step $1$.


\textbf{Step 2.}  At each crossing, there are two families of loops
without branch points and with exactly one corner.  Orient one of them
by $\left< \partial_\tau \right>$. As we shall see below, it does not
matter which, though for definiteness we assume that it is the family
that goes around the quadrant about which the orientation of the knot
runs counter-clockwise.

\textbf{Step 3.}  The third and final step in constructing a coherent
system of orientations relies on a diagrammatic understanding of the
gluing process and its relationship with orientations.  Let $\Sigma$
and $\Sigma'$ be families of loops that have corners at a common
crossing.  Suppose that $\Sigma$ is oriented by $\left< v_1, \ldots,
  v_m, \partial_\tau \right>$ and that $\Sigma'$ is oriented by
$\left< v_1', \ldots, v_n', \partial_\tau' \right>$.  Suppose further
that loops in $\Sigma$ go around a positive quadrant and loops in
$\Sigma'$ go around a negative quadrant.\footnote{In other words,
  $\Sigma$ lies ``below'' $\Sigma'$ in the symplectization picture.}
In order to understand the orientation on the glued family of loops,
we will examine what happens to the original families when we glue
them together geometrically near the corners.

\begin{figure}
  \centerline{\includegraphics{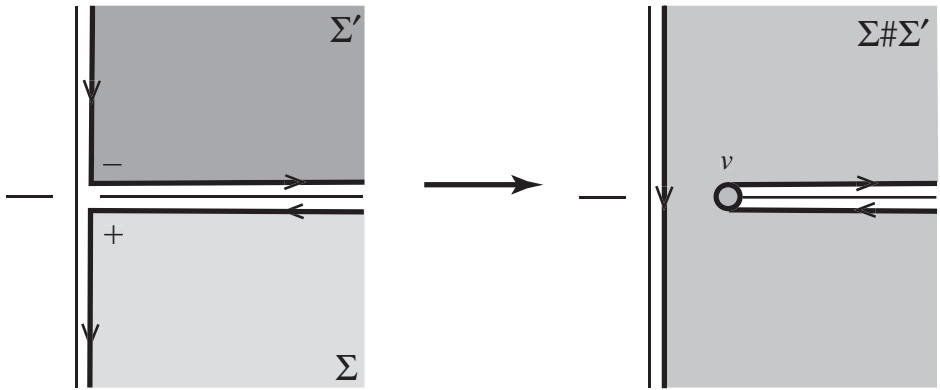}}
  \caption{The gluing process as viewed in the $xy$ projection.  Note
    that $\alpha$ parametrizes the position of the branch point in the
    right-hand diagram.}  \mylabel{fig:glue}
\end{figure}

Diagrammatically, gluing $\Sigma$ to $\Sigma'$ creates a new branch
point, as shown in Figure~\ref{fig:glue}.  We want to make this
process analytically precise in the ``straight line'' model that we
first examined in the proof of the Translation Theorem
\ref{thm:translation}.  In fact, it suffices to consider the straight
line model since, after a homotopy, we may assume that the operators
represented by the loops in $\Sigma$ and $\Sigma'$ come from the
straight line model near the crossing. In this model, $\Sigma$ is
parametrized by the vertical translation factor $a \in \rr$ and
consists of the maps
\begin{equation}
  f_a (u,v): (u,v) \mapsto \left(e^{-u}\cos(v), -e^{-u}\sin(v), v,
    u+\frac{e^{-2u}}{2}\cos^2(v)+a \right).
\end{equation}
Similarly, $\Sigma'$ is parametrized by $a' \in \rr$ and consists of
the maps
\begin{equation}
  f_{a'} (u,v): (u,v) \mapsto \left(e^{u}\cos(v), e^{u}\sin(v), v,
    u+\frac{e^{2u}}{2}\cos^2(v)+a' \right).
\end{equation}

We claim that the family of loops given locally by
\begin{multline}
  F_{\alpha, \tilde{a}}: (u,v) \mapsto \left(e^{-\alpha}\cosh(u)
    \cos(v), e^{-\alpha} \sinh(u) \sin(v), v, \right. \\
    \left. u+\frac{e^{-2\alpha}}{2} \cosh^2(u) \cos^2(v) + \tilde{a} \right)
\end{multline}
is the result of gluing $\Sigma$ to $\Sigma'$.  This family is
parametrized by the pair $(\alpha, \tilde{a})$.  The inverse of the
gluing map is given locally by:
\begin{equation}
  \Psi: (\alpha, \tilde{a}) \mapsto (\alpha + \tilde{a}, -\alpha +
  \tilde{a}).
\end{equation}
The important point is that $\Psi$ is orientation
preserving, so the orientation on the glued family given by
\begin{equation}\mylabel{eqn:glued-ori}
  \sigma(\Sigma) \# \sigma(\Sigma') =
  \left<v_1, \ldots, v_m, v_\alpha, v_1', \ldots, v_n', \partial_{\tilde{\tau}}
  \right>
\end{equation}
is the correct one.

We now proceed to orient all remaining families of loops using the
gluing procedure and the orientations chosen in steps $1$ and $2$ (and
on the trivial strip).  The idea is to start with a family of loops
with corners, then glue to the families from step $2$ until there are
no corners left, and finally compare the result to the orientations in
step $1$.

First, we orient all of the
families with one branch (or no) point(s) that lie around any given
crossing (see Figure~\ref{fig:vertex-loops}).  For now, assume that
the crossing in question is coherent about a $+$ (see
Section~\ref{ssec:alg}).  Glue $\Sigma_1$ to $\Sigma_2$, as pictured
in Figure~\ref{fig:glue-basic}.  If we orient $\Sigma_2$ with $\left<
  \pm w, \partial_\tau \right>$, then the glued loop is oriented by
$\left< v, \pm w, \partial_\tau \right>$.  Observe that $v$ does not
have the proper alignment to be listed first.  The orientation for the
corner-free family in step $1$ is $\left< w, v, \partial_\tau
\right>$, so, upon comparison, the orientation $\sigma(\Sigma_2)$ must
be $\left< -w, \partial_\tau \right>$.

Orient the family $\Sigma_4$ using the same procedure as above.  The
resulting orientation is $\left<w, \partial_\tau \right>$.

\begin{figure}
  \centerline{\includegraphics{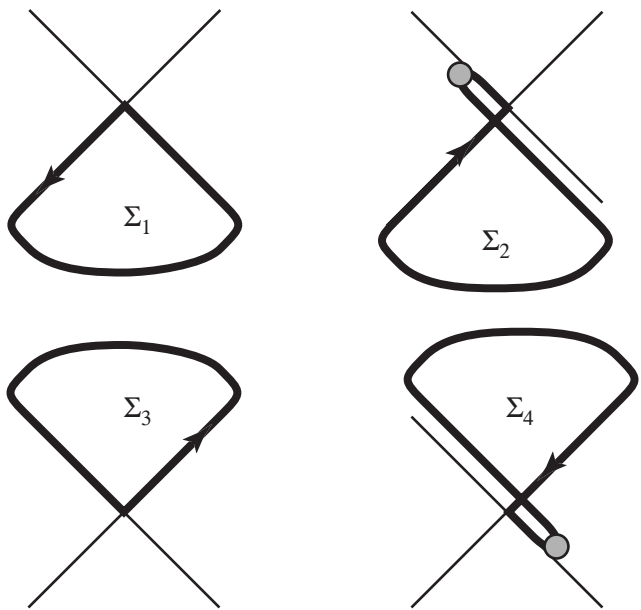}}
  \caption{A list of the four loops that lie around a vertex.  Each
    loop traverses half of $K$.}  \mylabel{fig:vertex-loops}
\end{figure}

\begin{figure}
  \centerline{\includegraphics{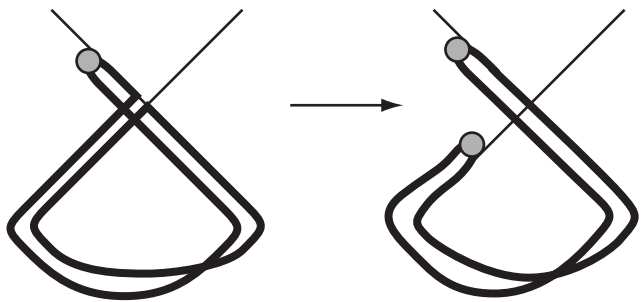}}
  \caption{Gluing a basic loop to a $(1+1)$-dimensional loop with a corner
    at a `$-$'.}  \mylabel{fig:glue-basic}
\end{figure}

Finally, to check that the choice of coherent quadrant in step $2$
does not matter, glue $\Sigma_2$ to the top family $\Sigma_3$, as
shown in Figure~\ref{fig:glue-basic2}.  If we orient $\Sigma_3$ by
$\left< \pm \partial_\tau \right>$, then the glued family is oriented
by $\left< \pm v, -w, \partial_\tau \right>$.  Since alignment
considerations show that the appropriate orientation on the
$(2+1)$-dimensional vertex-free family from step $1$ is $\left< w, v,
  \partial_\tau \right>$, the orientation on $\Sigma_3$ should be
$\left< \partial_\tau \right>$.  That this is the same orientation as
for $\Sigma_1$ comes from the fact that $\Sigma_1$ and $\Sigma_3$
represent the same family of operators.

So far, we have oriented all families with one corner at a vertex that
is coherent about a $+$.  To orient families around a vertex that is
coherent about a $-$, we follow exactly the same procedure, taking
care to use the ``bottom-to-top'' order of gluing.  As a result, the
signs on the vector $w$ for the loops covering the incoherent
quadrants are reversed: $\left< w, \partial_\tau \right>$ for the
leftmost (outward-pointing) family $\Sigma_2$ and $\left<-w,
  \partial_\tau \right>$ for the rightmost (inward-pointing) family
$\Sigma_4$. On the other hand, both $(0+1)$-dimensional families are
still oriented by $\left< \partial_\tau \right>$.

\begin{figure}
  \centerline{\includegraphics{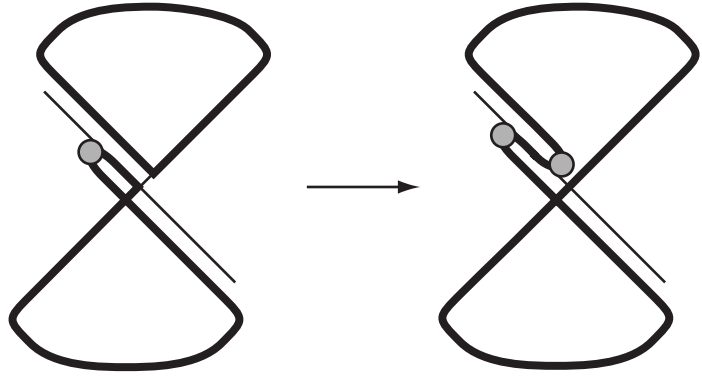}}
  \caption{Gluing the $(1+1)$-dimensional loop with a corner at a `$-$' to
    the loop opposite the basic loop.}  \mylabel{fig:glue-basic2}
\end{figure}

In the general case, we consider a family of loops $\Sigma$ with
positive corners labeled by $i_1, \ldots, i_m$ (in a cyclic order)
and negative corners labeled by $j_1, \ldots, j_n$.  Suppose that
$\Sigma$ has an orientation given by
\begin{equation}
  \sigma(\Sigma) = \left<v_1, \ldots, v_a, \pm\partial_\tau \right>.
\end{equation}
Our task is to determine the sign on $\partial_\tau$ that will fit
$\sigma(\Sigma)$ into a coherent system of orientations.  At each
crossing $i_j$, glue in one of the families $\Sigma_i$.  After gluing
to all of the corners, we are left with a family that has no corners
and an orientation given by:
$$\sigma(\Sigma_{i_m}) \# \cdots \# \sigma(\Sigma_{i_1}) \#
\sigma(\Sigma) \# \sigma(\Sigma_{j_1}) \# \cdots \#
\sigma(\Sigma_{j_n}).$$
Compare this orientation to the orientation
constructed in step $1$ and choose the sign of $\partial_\tau$ so that
the two orientations agree.  This completes step $3$.

\subsection{The Algorithm for Dimension $(0+1)$ Disks}
\mylabel{ssec:algorithm}

In this section, we explain how to obtain signs for dimension $(0+1)$
families of loops that represent immersed disks.  Our goal is to
translate the gluing process in the previous section into a concrete,
computable algorithm that, as expected, depends only on contributions
from the corners.  We get an actual sign rather than just an
orientation for the family $\Sigma(a; b_1, \ldots, b_k)$ by comparing
the orientation $\sigma(\Sigma)$ with the ``flow orientation'' $\left<
  \partial_\tau \right>$. Specifically, we show:
\begin{theorem}
        The sign on a holomorphic disk $f\in\ms^A/\rr$ in a moduli space of
        dimension $(0+1)$ given by comparing the coherent orientation on
        $f$ and the ``flow orientation'' is the same as the sign on $p(f)\in\Delta^A$
        described in Definition~\ref{def:signs}.
\end{theorem}

The remainder of this section is devoted to the proof of this theorem.

Our analysis begins with step 3 of the previous section.  The idea is
to look at the contributions that come from gluing to each corner
individually.  The order in which we glue depends on the
counter-clockwise cyclic order in the corners, i.e. we start with $a$
and end up with $b_k$.  At a coherent corner, gluing gives two new
vectors in the orientation (see Figure~\ref{fig:glue-basic2}, for
example); at an incoherent corner, gluing gives only one new vector.

At a coherent corner of either sign, the sign contribution depends on
the number of negative vectors added during the gluing process.  If
there is one negative vector, then, ignoring alignment for the moment,
we need to flip a single sign to make the pair agree with the base
orientation from step $1.$  Thus, there is a sign contribution of $-1$
in this case.  On the other hand, the same reasoning shows that there
is no sign contribution when the coherent corner contributes an even
number of negative vectors.  See Figure~\ref{fig:coherent} for the
vectors given by the conventions adopted in the previous section.  It
follows from Figure~\ref{fig:coherent} that corners that are clockwise
at a `$+$' and counter-clockwise at a `$-$' pick up a negative sign.

\begin{figure}
  \centerline{\includegraphics[width=3.2in]{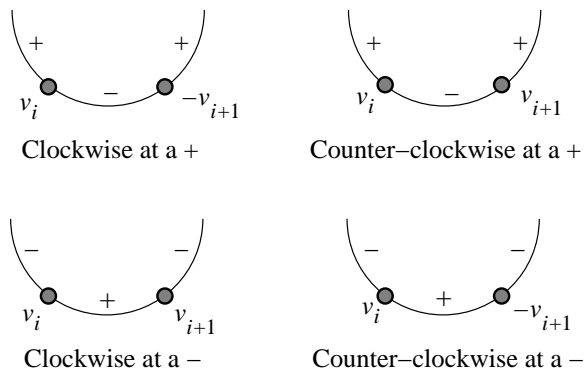}}
  \caption{Orientation choices at coherent vertices, pulled back to
    $S^1$ for clarity.}
  \mylabel{fig:coherent}
\end{figure}

Incoherent corners come in inward / outward pairs. Our analysis of
gluing shows that each pair contributes two positive vectors.  Hence,
by the same analysis as in the coherent case, there is no sign
contribution at an incoherent corner.

Finally, we consider alignment issues. Suppose that the positive
corner, $a$, is either clockwise coherent or inward-pointing.  On the
first segment of the loop, the orientation of the knot disagrees with
the counter-clockwise orientation of the family.  Since we agreed to
start labeling the branch points with one that changes alignment from
$+$ to $-$, we have to take into account a cyclic permutation of the
labelings.  This is an \emph{odd} permutation since a cornerless
family of loops has an even number of branch points.  Thus, we pick up
an extra negative sign at a $+$ vertex when it is either clockwise
coherent or inward-pointing.

In sum, we have deduced that Definition~\ref{def:signs} assigns the
correct sign to a dimension $(0+1)$ disk.

Note that there are actually sixteen different sets of choices for the
signs in Figure~\ref{fig:lenny-signs}.  This is reflected by the
orientations we chose in Section~\ref{ssec:geom-ori} on $D^2 \times
\cc$, $\mathcal{O}(1)$, the capping disk without branch points, and
the trivial vertical strip.  It turns out that all of the orientation
choices give the same DGA, possibly after replacing some generators by
their negatives.

\section{Acknowledgments}
\mylabel{sec:acks}

We would like to thank Yasha Eliashberg and Frederic Bourgeois for
many stimulating discussions about the material presented in this
paper.  In addition, we would like to thank the American Institute of
Mathematics (AIM) for facilitating the Low-Dimensional Contact
Geometry program, during which we completed this work.


\bibliography{oriented-chvII-rev}

\end{document}